\documentclass[a4paper]{article}
\usepackage{amssymb,amscd,amsfonts,mathrsfs,amsmath}
\usepackage[all]{xy}
\input amssym.def

\def\double{\mathbb}

\def\cc{{\double C}}

\def\zz{{\double Z}}

\def\rr{{\double R}}

\newtheorem{theorem}{Theorem}[section]
\newtheorem{lemma}[theorem]{Lemma}
\newtheorem{corollary}[theorem]{Corollary}
\newtheorem{definition}[theorem]{Definition}
\newtheorem{proposition}[theorem]{Proposition}

\newtheorem{example}[theorem]{Example}

\def\Uc{{\mathscr U}}
\def\res{\mathop{\mathrm{Res}}\limits_{z=0}}

\def\cp{\rtimes}
\def\si{\sigma}

\def\cinf{C^{\infty}}
\def\cinfc{C^{\infty}_c}

\newcommand{\be}{\begin{equation}}
\newcommand{\ee}{\end{equation}}
\newcommand{\beq}{\begin{eqnarray}}
\newcommand{\eeq}{\end{eqnarray}}
\newcommand{\om}{\omega}
\newcommand{\Om}{\Omega}
\newcommand{\al}{\alpha}
\def\nat{\natural}
\def\id{{\mathop{\mathrm{id}}}}

\newcommand{\La}{\Lambda}
\newcommand{\la}{\lambda}
\newcommand{\Ec}{{\mathscr E}}
\newcommand{\Vc}{{\mathscr V}}
\newcommand{\Lc}{{\mathscr L}}
\newcommand{\non}{\nonumber}
\newcommand{\eps}{\varepsilon}
\newcommand{\Sc}{{\mathscr S}}
\newcommand{\Wc}{{\mathscr W}}
\newcommand{\Yc}{{\mathscr Y}}
\newcommand{\Zc}{{\mathscr Z}}
\newcommand{\Rc}{{\mathscr R}}

\newcommand{\Ic}{{\mathscr I}}
\newcommand{\Jc}{{\mathscr J}}
\newcommand{\Oc}{{\mathscr O}}
\newcommand{\Ind}{{\mathop{\mathrm{Ind}}}}

\def\Td{\mathrm{Td}}
\def\ad{\mathrm{Ad\,}}

\newcommand{\Tr}{{\mathop{\mathrm{Tr}}}}

\newcommand{\tr}{{\mathop{\mathrm{tr}}}}
\newcommand{\Ac}{{\mathscr A}}
\newcommand{\Gc}{{\mathscr G}}
\newcommand{\te}{\theta}

\newcommand{\Te}{\Theta}

\newcommand{\cqfd}{\hfill\rule{1ex}{1ex}}

\def\Id{\mathrm{Id}}

\def\d{\partial}
\def\dd{\mathrm{\bf d}}
\def\DD{\mathrm{\bf D}}

\def\Hc{{\mathscr H}}

\def\Bc{{\mathscr B}}
\def\Cc{{\mathscr C}}
\def\Jc{{\mathscr J}}

\def\Fc{{\mathscr F}}
\def\Pc{{\mathscr P}}
\def\Qc{{\mathscr Q}}
\def\im{\mathop{\mathrm{Im}}}
\def\ker{\mathop{\mathrm{Ker}}}

\def\bb{\overline{b}}
\def\phib{\overline{\phi}}
\def\psib{\overline{\psi}}

\def\hom{{\mathop{\mathrm{Hom}}}}
\def\dom{{\mathop{\mathrm{Dom}}}}

\def\End{{\mathop{\mathrm{End}}}}

\def\Tt{\widetilde{T}}

\def\et{\tilde{e}}

\def\Omh{\widehat{\Omega}}
\def\Th{\widehat{T}}

\def\Pch{\widehat{\mathscr P}}
\def\Hch{\widehat{\mathscr H}}
\def\Vch{\widehat{\mathscr V}}
\def\Ych{\widehat{\mathscr Y}}

\def\Xc{{\mathscr X}}

\def\Tc{{\mathscr T}}
\def\Dc{{\mathscr D}}
\def\Uc{{\mathscr U}}

\def\supp{\mathrm{supp}\,}

\def\uh{\hat{u}}

\def\ah{\hat{a}}

\def\sign{\mathrm{sign}}
\def\Phib{\overline{\Phi}}
\def\nablab{\overline{\nabla}}
\def\CL{\mathrm{CL}}
\def\LS{\mathrm{LS}}

\def\CS{\mathrm{CS}}
\def\PS{\mathrm{PS}}

\def\i{\mathrm{i}}

\def\R{\mathrm{R}}
\def\top{\mathrm{top}}

\def\Res{\mathrm{Res}}

\def\Ad{\mathrm{Ad}}

\begin{document}

\begin{center}

{\bf INDEX THEORY FOR IMPROPER ACTIONS:\\ LOCALIZATION AT UNITS}
\vskip 1cm
{\bf Denis PERROT}
\vskip 0.5cm
Institut Camille Jordan, CNRS UMR 5208\\
Universit\'e de Lyon, Universit\'e Lyon 1,\\
43, bd du 11 novembre 1918, 69622 Villeurbanne Cedex, France \\[2mm]
{\tt perrot@math.univ-lyon1.fr}\\[2mm]
\end{center}
\vskip 0.5cm
\begin{abstract}
We pursue the study of local index theory for operators of Fourier-integral type associated to non-proper and non-isometric actions of Lie groupoids, initiated in a previous work. We introduce the notion of geometric cocycles for Lie groupoids, which allow to represent fairly general cyclic cohomology classes of the convolution algebra of Lie groupoids localized at isotropic submanifolds. Then we compute the image of geometric cocycles localized at units under the excision map of the fundamental pseudodifferential extension. As an illustrative example, we prove an equivariant longitudinal index theorem for a codimension one foliation endowed with a transverse action of the group of real numbers. 
\end{abstract}

\vskip 0.5cm


\section{Introduction}

In a previous paper \cite{P10} we initiated a local index theory for Fourier-integral type operators associated to non-proper and non-isometric actions of Lie groupoids on smooth submersions. The difficulty in the non-isometric setting is that the main tools of ordinary local index theory, such as the heat kernel expansion or related methods, are no longer available. One thus has to invent entirely novel techniques in order to handle such general situations.\\
The precise setting is the following. Let $G\rightrightarrows B$ be a (Hausdorff) Lie groupoid acting smoothly on a submersion of manifolds $\pi:M\to B$. We denote by $\cinfc(B,\CL_c(M))$ the algebra of compactly-supported classical pseudodifferential operators acting along the fibers of the submersion. This algebra naturally carries an action of $G$, so we can form the crossed product $\cinfc(B,\CL_{c}(M))\cp G$. The latter contains many interesting operators, including of Fourier-integral type, that are of interest from the point of view of non-commutative index theory. We can form the short exact sequence of algebras 
$$
0 \to \cinfc(B,\CL^{-1}_c(M))\cp G \to \cinfc(B,\CL^0_c(M))\cp G \to \cinfc(S^*_{\pi}M)\cp G \to 0
$$
(see section \ref{sres}). The quotient $\cinfc(S^*_{\pi}M)\cp G$ is isomorphic to the convolution algebra of the action groupoid $S^*_{\pi}M\cp G$, and describes the non-commutative leading symbols of the operators under consideration. The fundamental question is to compute the excision map $E^*$ induced by this extension on periodic cyclic cohomology. We solve this problem in \cite{P10} in the case of cyclic cohomology classes which are \emph{localized} at appropriate isotropic submanifolds $O\subset G$. The main result of \cite{P10}, which we shall recall here in section \ref{sres} (see Theorem \ref{tres}),  takes the form of a commutative diagram:
\be 
\vcenter{\xymatrix{
HP^\bullet(\cinfc(B,\CL_{c}^{-1}(M))\cp G) \ar[r]^{\quad E^*} & HP^{\bullet+1}(\cinfc(S^*_\pi M)\rtimes G)  \\
HP_\top^{\bullet}(\cinfc(B)\cp G)_{[O]} \ar[u]^{\tau_*} \ar[r]^{\pi^!_G\qquad} & HP_\top^{\bullet+1}(\cinfc(S^*_\pi M)\rtimes G)_{[\pi^* O]} \ar[u] }}
\ee
The bottom map $\pi^!_G$ is given by a fixed-point formula, involving residues of zeta-functions which generalize the well-known Wodzicki residue \cite{Wo}.\\

Although the residue formula is an explicit local expression in the complete symbols of the operators, it is extremely hard to compute. The reason is that when the dimension $n$ of the fibers of the submersion $\pi: M\to B$ is ``large" (typically $n \geq 2$ !) the residues involve an asymptotic expansion of symbols up to order $n$, which produces a huge quantity of terms. However, in the case of cyclic cohomology \emph{localized at units} ($O=B$), using the techniques of \cite{P9} we shall give here a close expression for the map $\pi^!_G$ in terms of the usual characteristic classes entering the index theorem, namely the $G$-equivariant Todd class of the vertical tangent bundle of $M$. To show this we focus on a class of elements in $HP^\bullet_\top(\cinfc(B)\cp G)_{[O]}$ represented by \emph{geometric cocycles}: the latter are quadruples $(N,E,\Phi,c)$ where $N,E$ are smooth manifolds, $\Phi$ is a flat connection on a certain groupoid, and $c$ is a cocycle in a certain complex. See section \ref{sgeo} for the precise definitions. This geometric construction of localized cylic cohomology classes is of independent interest since it works for the cyclic cohomology localized at any isotropic submanifold $O\subset G$. Moreover in the case of the localization at units $O=B$, one recovers all the previously known constructions of cyclic cohomology in special cases, including: the cohomology of the classifying space for $G$ and Gelfand-Fuchs cohomology (\'etale groupoids or foliation groupoids \cite{C83, C94}), differentiable groupoid cohomology (\cite{PPT}), etc. The next result, Theorem \ref{tgeo} of this article, is a refinement of Theorem \ref{tres} in the case of geometric cocycles localized at units.\\

\noindent {\bf Theorem 5.4} {\it Let $G\rightrightarrows B$ be a Lie groupoid acting on a surjective submersion $\pi:M\to B$. The excision map localized at units
$$
\pi^!_G\ :\ HP_\top^{\bullet}(\cinfc(B)\cp G)_{[B]} \to HP_\top^{\bullet+1}(\cinfc(S^*_\pi M\rtimes G))_{[S^*_\pi M]}
$$
sends the cyclic cohomology class of a \emph{proper} geometric cocycle $(N,E,\Phi,c)$ to the cyclic cohomology class
$$
\pi^!_G([N,E,\Phi,c]) = [N\times_B S^*_\pi M \, , \, E\times_B S^*_\pi M \, , \,  \pi_*^{-1}(\Phi) \, , \, \Td(T_\pi M\otimes\cc) \wedge \pi^*(c)] 
$$
where $\Td(T_\pi M\otimes\cc)$ is the Todd class of the complexified vertical tangent bundle in the invariant leafwise cohomology of $E\times_B S^*_\pi M$.}\\

\noindent Once specialized to the case of a discrete group $G$ acting by \emph{any} diffeomorphisms on a manifold $M$, this theorem completely solves the problem of evaluating the $K$-theoretical index of the Fourier integral operators considered in \cite{SS} on cyclic cohomology classes localized at units, in terms of the non-commutative leading symbol of this operator in $K_1(\cinf(S^*M)\cp G)$. Let us mention that an adaptation of this result to the hypoelliptic calculus of Connes and Moscovici (see \cite{CM95}) gives a direct proof of their index theorem for \emph{transversally} hypoelliptic operators on foliations, and moreover gives an explicit computation of the characteristic class of the hypoelliptic signature operator. This is the topic of a separate paper \cite{PR}. \\
To end this article we also present an application of the residue formula of Theorem \ref{tres} in a case \emph{not} localized at units. We prove an equivariant longitudinal index theorem for a codimension one foliation endowed with a transverse action of the group $\rr$. Choosing a complete transversal for the foliation, this geometric situation can be reduced to the action of a Lie groupoid on a submersion. The groupoid possesses a canonical trace, and we show that the pairing of the corresponding cyclic cocycle with the index of any leafwise elliptic equivariant pseudodifferential operator is given by a formula localized at the periodic orbits of the transverse flow on the foliation. This gives an interesting interpretation of the results of Alvarez-Lopez, Kordyukov \cite{ALK} and Lazarov \cite{L} in the context of the $K$-theory/cyclic cohomology of Lie groupoids.\\

Here is a brief description of the article. Section \ref{sres} recalls the main result of \cite{P10}, namely the residue theorem computing the excision map in cyclic cohomology associated to the fundamental pseudodifferential extension. Section \ref{sgeo} describes a very general construction of cyclic cohomology classes from geometric cocycles. In the particular case of localization at units, this construction is related to the computation of the excision map in section \ref{sloc}. Then Theorem \ref{tgeo} is proved in section \ref{sdir}, adapting the abstract Dirac construction of \cite{P9}. Finally section \ref{sfol} contains the example of a transverse flow on a codimension one foliation. \\

\section{The residue theorem}\label{sres}

Let $G\rightrightarrows B$ be a (Hausdorff) Lie groupoid, with $B$ its manifold of units. We denote by $r:G\to B$ the range map and by $s:G\to B$ the source map. Both are smooth submersions. We think of an element $g\in G$ as a left-oriented arrow:
$$
g: r(g) \longleftarrow s(g)
$$
The restriction of the vector bundle $\ker s_*\subset TG$ to the submanifold $B\subset G$ is the Lie algebroid $AG$ over $B$. Denote by $A^*G$ its dual vector bundle. Under the range map, the line bundle $\vert\Lambda^{\max} A^*G \vert $ over $B$ can be pulled-back to a line bundle $r^*(\vert\Lambda^{\max} A^*G \vert)$ over $G$. The latter is canonically identified with the bundle of 1-densities along the fibers of the source map. By definition, the convolution algebra of $G$ is the space of smooth compactly supported sections of this line bundle (complexified):
\be
\cinfc(B)\cp G := \cinfc(G, r^*(\vert\Lambda^{\max} A^*G \vert))\ .
\ee
The product of two sections $a_1$ and $a_2$ is given by convolution,
\be
(a_1a_2)(g) = \int_{g_1g_2=g} a_1(g_1)a_2(g_2)
\ee
where the integral is taken with respect to the 1-densities. More generally, if $\R$ is a $G$-equivariant vector bundle over $B$, we denote by $U_g: \R_{s(g)}\to \R_{r(g)}$ the vector space isomorphism induced by an element $g\in G$ between the fibers of $\R$ above the points $s(g)$ and $r(g)$. If moreover $\R$ is a $G$-equivariant algebra bundle, we define as above the convolution algebra 
\be
\cinfc(B,\R)\cp G := \cinfc(G, r^*(\R\otimes \vert\Lambda^{\max} A^*G \vert))\ ,
\ee
where the convolution product is now twisted by the action of $G$ on $\R$: 
\be
(a_1a_2)(g) = \int_{g_1g_2=g} a_1(g_1) \cdot U_{g_1}a_2(g_2)
\ee
Let now the groupoid $G\rightrightarrows B$ act (from the right) on a smooth submersion of manifolds $\pi: M\to B$. In \cite{P10} we define the $G$-equivariant algebra bundle $\CL^0_c(M)$, whose fiber over a point $b\in B$ is the algebra of compactly supported \emph{classical} (one-step polyhomogeneous) pseudodifferential operators of order $\leq 0$ on the submanifold $M_b=\pi^{-1}(b)$. It contains the algebra subbundle $\CL^{-1}_c(M)$ of order $\leq -1$ operators as a two-sided ideal, whence the fundamental extension $(E)$ of convolution algebras
$$
0 \to \cinfc(B,\CL^{-1}_c(M))\cp G \to \cinfc(B,\CL^0_c(M))\cp G \to \cinfc(S^*_{\pi}M)\cp G \to 0
$$
where $S^*_{\pi}M$ is the vertical sphere bundle over $M$, and the epimorphism onto the convolution algebra $\cinfc(S^*_{\pi}M)\cp G$ is induced by the leading symbol. The goal is to calculate the connecting map induced by this extension on periodic cyclic cohomology
\be
E^*\ :\  HP^\bullet(\cinfc(B,\CL_{c}^{-1}(M))\cp G) \to HP^{\bullet+1}(\cinfc(S^*_\pi M)\rtimes G)\ .
\ee
As it turns out, we can give explicit residue formulas for the map $E^*$ when cyclic cohomology is localized. An isotropic submanifold $O\subset G$ is a submanifold with the property that any element $g\in O$ verifies $r(g)=s(g)$. It is $\Ad$-invariant if it is stable by conjugation by any element of $G$. In  \cite{P10} we define the (topological) cyclic cohomology of $G$ localized at any $\Ad$-invariant isotropic submanifold $O$, which comes with a natural forgetful map to ordinary (algebraic) cyclic cohomology
$$
HP_\top^{\bullet}(\cinfc(B)\cp G)_{[O]}  \to HP^{\bullet}(\cinfc(B)\cp G)\ .
$$
Note that an element $g\in O$ acts by diffeomorphism on the manifold $M_{r(g)}=\pi^{-1}(r(g))$. We say that the action of $O$ on the submersion $M$ is \emph{non degenerate} if 
\begin{itemize}
\item For any $g\in O$, the set of fixed points $M^g_{r(g)}$ is a union of isolated submanifolds in $M_{r(g)}$, depending smoothly on $g$;
\item At any point $x\in M_{r(g)}^g$ the tangent space $T_xM_{r(g)}$ in the ambient manifold $M_{r(g)}$ splits as a direct sum
\be
T_xM_{r(g)} = T_xM_{r(g)}^h\oplus N_x^g
\ee
of two subspaces globally invariant by the action of the tangent map $g_*$ associated to the diffeomorphism. We denote $g'$ the restriction of $g_*$ to the normal subspace $N_x^g$;
\item The endomorphism $1-g'$ of $N_x^g$ is non-singular, that is $\det(1-g')\neq 0$ at any point $x\in M_{r(g)}^g$.
\end{itemize}
We then have the following residue theorem:
\begin{theorem}[\cite{P10}]\label{tres}
Let $G\rightrightarrows B$ be a Lie groupoid and let $O$ be an $\Ad$-invariant isotropic submanifold of $G$. Let $\pi:M\to B$ be a $G$-equivariant surjective submersion and assume the action of $O$ on $M$ non-degenerate. Then one has a commutative diagram 
\be 
\vcenter{\xymatrix{
HP^\bullet(\cinfc(B,\CL_{c}^{-1}(M))\cp G) \ar[r]^{\quad E^*} & HP^{\bullet+1}(\cinfc(S^*_\pi M)\rtimes G)  \\
HP_\top^{\bullet}(\cinfc(B)\cp G)_{[O]} \ar[u]^{\tau_*} \ar[r]^{\pi^!_G\qquad} & HP_\top^{\bullet+1}(\cinfc(S^*_\pi M)\rtimes G)_{[\pi^* O]} \ar[u] }}
\ee
where the isotropic submanifold $\pi^*O \subset S^*_\pi M\rtimes G$ is the pullback of $O$ by the submersion $S^*_\pi M\to B$. The map $\pi^!_G$ is expressed by an explicit residue formula.
\end{theorem}
The left vertical arrow comes from the canonical Morita equivalence between the groupoid $G$ and its pullback $\pi^*G$ under the submersion $\pi$. \\

Let us now explain the way $\pi^!_G$ is constructed. For convenience we write $\Ac = \cinfc(S^*_\pi M)\cp G$, $\Ec = \cinfc(B,\CL_c^0(M))\cp G$, $\Bc = \cinfc(B,\CL_c^{-1}(M))\cp G$. We use the Cuntz-Quillen formalism \cite{CQ95} and compute the cyclic homology of $\Ac$ by means of the $X$-complex of its non-unital tensor algebra $T\Ac$, completed in the adic topology relative to the ideal $J\Ac=\ker(T\Ac\to\Ac)$:
$$
\Th\Ac = \varprojlim_n T\Ac /J\Ac^n \ ,\qquad X(\Th\Ac) : \Th\Ac \rightleftarrows \Om^1\Th\Ac_\nat \ .
$$
Then $HP^\bullet(\Ac)$ is the cohomology of the $\zz_2$-graded complex $\hom(X(\Th\Ac),\cc)$. The extension $0 \to \Bc \to \Ec \to \Ac \to 0$ can be inserted in a commutative diagram where all rows and columns are exact:
\be
\vcenter{\xymatrix{
 & 0 \ar[d] & 0 \ar[d] & 0 \ar[d] &  \\
0 \ar [r] & \Jc \ar[r] \ar[d] & J\Ec \ar[r] \ar[d] & J\Ac  \ar[r] \ar[d] & 0  \\
0 \ar [r] & \Rc \ar[r] \ar[d] & T\Ec \ar[r] \ar[d] & T\Ac  \ar[r] \ar[d] & 0  \\
0 \ar [r] & \Bc \ar[r] \ar[d] & \Ec \ar[r] \ar[d] & \Ac  \ar[r] \ar[d] & 0  \\
 & 0 & 0 & 0 &  }}\label{diag}
\ee
The ideal $\Rc$ is the kernel of the homomorphism $T\Ec\to T\Ac$, and $\Jc$ is the kernel of $J\Ec\to J\Bc$. The kernel of the homomorphism $T\Ec\to\Ac$ is the ideal $J\Ec + \Rc$. We define the corresponding adic completion of $T\Ec$ and its $X$-complex:
$$
\widetilde{T}\Ec = \varprojlim_n T\Ec/(J\Ec+\Rc)^n\ ,\qquad X(\widetilde{T}\Ec) : \Tt\Ec \rightleftarrows \Om^1\Tt\Ec_\nat \ .
$$
Choose a continuous linear splitting $\si:\Ac\to \Ec$ of the leading symbol homomorphism. As a vector space, $\Ec$ is a direct summand in $T\Ec$, hence we can view $\si$ as a linear map to $T\Ec$. By the universal property of the tensor algebra $T\Ac$ we can lift $\si$ to a homomorphism of algebras $\si_*:T\Ac\to T\Ec$:
$$
\vcenter{\xymatrix{0 \ar[r] & J\Ac \ar[r] \ar[d]_{\si_*} & T\Ac \ar[r] \ar[d]_{\si_*} & \Ac \ar[r] \ar@{.>}[dl]^{\si} \ar@{=}[d] & 0 \\
0\ar[r] & J\Ec + \Rc \ar[r]  & T\Ec \ar[r] & \Ac \ar[r]   & 0 }}
$$
One has $\si_*(a_1\otimes\ldots\otimes a_n) = \si(a_1)\otimes \ldots \otimes \si(a_n)$ in $T\Ec$. Because $\si_*$ respects the ideals, it extends to a homomorphism $\Th\Ac\to\widetilde{T}\Ec$ and also induces a chain map 
$$
\si_*:X(\Th\Ac)\to X(\widetilde{T}\Ec)\ .
$$
Let $\Cc=\cinfc(B)\cp G$ be the convolution algebra of $G$, and $\Th\Cc$ be the $J\Cc$-adic completion of its tensor algebra. Taking the locally convex topology of $\Cc$ into account, we can replace algebraic tensor products everywhere by topological ones and get the larger algebra $\Th\Cc_\top$. One can think of an element in $\Th\Cc_\top$ as a collection of smooth compactly supported functions over the product manifolds $G^n$ for all $n$. The localization $X(\Th\Cc_\top)_{[O]}$ of its $X$-complex at an $\Ad$-invariant isotropic submanifold $O\in G$ is defined in \cite{P10} as a quotient, and its elements may be viewed as collections of jets of smooth functions at localized submanifolds in $G^n$. The localized cyclic cohomology $HP_\top^{\bullet}(\cinfc(B)\cp G)_{[O]}$ is the cohomology of the $\zz_2$-graded complex $\hom(X(\Th\Cc_\top),\cc)$ of cochains which satisfy certain conditions of \emph{continuity} and \emph{boundedness}. There exists a \emph{residue morphism}
\be 
\Res\ :\ X(\Tt\Ec) \to X(\Th\Cc_\top)_{[O]} \label{resmap} \ .
\ee
This requires to choose a smooth section $Q\in \cinf(B,\CL^1(M))$ of elliptic, positive pseudodifferential operators acting along the fibers of the submersion $\pi$. Thus at any point $b\in B$ one has an elliptic positive invertible operator $Q_b$ acting on the manifold $M_b=\pi^{-1}(b)$. We may furthermore assume that, modulo perturbation by regularizing operators, each $Q_b$ as well as its complex powers $Q_b^{-z}$, $z\in \cc$ are \emph{properly supported}. For Re$(z)$ sufficiently large, $Q_b^{-z}$ is used to regularize operator traces and leads to meromorphic zeta-functions. In even degree (\ref{resmap}) is a linear map $\Tt\Ec\to (\Th\Cc_\top)_{[O]}$. The image of a tensor $e_1\otimes\ldots \otimes e_n \in T\Ec$ is a jet of a function of $n$ variables $(g_1,\ldots,g_n)$ in $G^n$: 
\beq
\lefteqn{\big(\Res(e_1\otimes\ldots \otimes e_n)\big)(g_1,\ldots,g_n) = } \non\\
&&\res \Tr\big(e_1(g_1)\, U_{g_1}\, {h(s(g_1),r(g_2))}  \ldots e_n(g_n)\, U_{g_n}\, {h(s(g_n),r(g_1))}\, Q^{-z}_{r(g_1)}\big) \non
\eeq
where $h$ is a \emph{generalized connection} on the submersion $\pi:M\to B$ according to Definition 4.3 of \cite{P10}. In odd degree (\ref{resmap}) is a linear map $\Om^1\Tt\Ec_\nat\to (\Om^1\Th\Cc_\top)_{[O]}$. It is defined in a similar way:
\beq
\lefteqn{\big(\Res(\nat(e_1\otimes\ldots \otimes e_{n-1}\dd e_n))\big)(g_1,\ldots,g_{n-1}|g_n) = } \non\\
&& \res \Tr\big(e_1(g_1)\, U_{g_1}\, {h(s(g_1),r(g_2))}  \ldots  Q^{-z}_{r(g_n)}\,  e_n(g_n)\, U_{g_n}\, {h(s(g_n),r(g_1))}\big) \non
\eeq
These residues select the poles of zeta-functions like $\Tr(PU_gQ_b^{-z})$, where $g\in O$, $b=s(g)=r(g)$, and $P$ is some classical pseudodifferential operator on the manifold $M_b=\pi^{-1}(b)$. To see that they yield in fact local formulas, Proposition 5.3 of \cite{P10} computes these residues in terms of integrals of the complete symbol $\si_P$ of $P$ over the cosphere bundle of the fixed-point set of $g$. Suppose for simplicity that $M_b^g\subset M_b$ is a unique submanifold of dimension $r$. Then there exists a local coordinate system $(x;y)$ on $M_b$ which is adapted to $M_b^g$ (see \cite{P10} for precise definitions), and a canonical coordinate system $(x,p;y,q)$ on the cotangent bundle $T^*M_b$, such that $(x,p)$ is a canonical coordinate system of $T^*M_b^g$. Then one has 
\be
\res \Tr(PU_g Q_b^{-z}) = \int_{S^*M_b^g} \left[ e^{\i\langle \frac{\d}{\d q}, (1-g')^{-1}\frac{\d}{\d y}\rangle} \cdot  \frac{ \si_P\, e^{\i\langle p, x-g^*x\rangle} }{|\det(1-g')|}  \right]_{-r} \, \frac{\eta (d\eta)^{r-1}}{(2\pi)^r} \label{floc}
\ee
where $S^*M_b^g$ is the cosphere bundle of the fixed submanifold, $\eta=\langle p , dx\rangle$ is the contact one-form over the cotangent bundle $T^*M_b^g$, the brackets $[\ ]_{-r}$ select the order $-r$ component of a symbol with respect to the variables $(x,p)$, and $g'$ is the matrix of partial derivatives $\d (g^*y)/\d y$.  In the special case $U_g=\Id$, this expression simplifies drastically and one recovers the well-known Wodzicki residue \cite{Wo}. Also note that in the general case the residue does not depend on the choice of $Q$. Now we enlarge $\CL_c(M)$ to the algebra bundle $\CL_c(M)_{\log}$ of \emph{log-polyhomogeneous} pseudodifferential operators and define as before the crossed product algebra
$$
\Ec_{\log} = \cinfc(B,\CL^0_c(M)_{\log})\cp G\ .
$$
$\Ec_{\log}$ contains products of classical operators in $\Ec$ by logarithms of elliptic pseudodifferential operators. For example the logarithm $\ln Q = - \frac{d}{dz}Q^{-z}|_{z=0}$ of $Q$ acts on $\Ec_{\log}$ as a left multiplier: $(\ln Q\cdot e)(g) = \ln Q_{r(g)} e(g)$ for all $e\in \Ec_{\log}$ and $g\in G$. Replacing all preceding algebras of pseudodifferential operators by log-polyhomogeneous ones, we define $\Tt\Ec_{\log}$ as the adic completion of $T\Ec_{\log}$ with respect to the ideal $J\Ec_{\log}+\Rc_{\log}$, and $X(\Tt\Ec)_{\log}^1$ as the  subcomplex of $X(\Tt\Ec)_{\log}$ having logarithmic degree at most $1$. Then, the residue map (\ref{resmap}) can be extended to a chain map
$$
X(\Tt\Ec)^1_{\log}\cap\dom(\Res) \to X(\Th\Cc_\top)_{[O]}
$$
with an appropriate domain $\dom(\Res)$. Typically, ``commutators" with $\ln Q$ in a tensor product belong to this domain, and will be written as
\beq
\lefteqn{e_1\otimes\ldots \otimes[\ln Q, e_i]\otimes e_{i+1}\otimes  \ldots \otimes e_n := }\non\\
&&  e_1\otimes\ldots \otimes \ln Q \cdot e_i \otimes e_{i+1}\otimes  \ldots \otimes e_n - e_1\otimes\ldots \otimes e_i \otimes \ln Q\cdot e_{i+1}\otimes  \ldots \otimes e_n \non
\eeq
for all $e_i\in \Ec$. Now everything is set to give the explicit formulas for the map $\pi^!_G$. The image of a localized cyclic cohomology class of even degree $[\varphi]\in HP_\top^{0}(\cinfc(B)\cp G)_{[O]}$ is represented by the following localized cyclic cocycle of odd degree $\pi^!_G(\varphi)\in \hom(\Om^1\Th\Ac_\nat,\cc)$:
$$
\pi^!_G(\varphi)(\nat (a_1\otimes \ldots\otimes a_{n-1} \dd a_n))  = \varphi\circ\Res\big( \si(a_1)\otimes \ldots \otimes\si(a_{n-1})\otimes [\ln Q, \si(a_n)] \big) \non
$$
for all $\nat(a_1\otimes \ldots\otimes a_{n-1}\dd a_n) \in \Om^1 T\Ac_\nat$. In a similar way, the image of a localized cyclic cohomology class of odd degree $[\varphi]\in HP_\top^{1}(\cinfc(B)\cp G)_{[O]}$ is represented by the following localized cyclic cocycle of even degree $\pi^!_G(\varphi)\in \hom(\Th\Ac,\cc)$:
$$
\pi^!_G(\varphi)(a_1\otimes \ldots\otimes a_n) = \!\!\!  \sum_{1\leq i<j\leq n} \!\!\! \varphi\circ\Res\big( \si(a_1)\otimes \ldots [\ln Q, \si(a_i)]\ldots \dd\si(a_j)\ldots \otimes\si(a_n)\big) 
$$
for all $a_1\otimes \ldots\otimes a_n \in T\Ac$. \\

By \cite{P10} Corollary 5.7, the residue theorem allows to evaluate the image of the index map in algebraic $K$-theory associated to the extension $E: 0\to \Bc \to \Ec \to \Ac \to 0$
\be
\Ind_E : K_1(\Ac) \to K_0(\Bc)
\ee
on localized cyclic cohomology classes. Let $[u]\in K_1(\Ac)$ be an elliptic symbol represented by an invertible matrix $u\in M_{\infty}(\Ac)^+$ in the unitalization of the matrix algebra $M_{\infty}(\Ac)$. It canonically lifts to an invertible element $\uh\in M_{\infty}(\Th\Ac)^+$ under the linear embedding $\Ac\subset \Th\Ac$. We extend the linear splitting $\si:\Ac\to\Ec$ to the unitalized matrix algebras $M_\infty(\Ac)^+$ and $M_\infty(\Ec)^+$ by setting $\si(1)=1$, and lift it to a unitalized algebra homomorphism $\si_*:M_\infty(\Th\Ac)^+\to M_\infty(\Tt\Ec)^+$. The latter maps $\uh$ to an invertible element $\si_*(\uh)$, whose inverse is 
$$
\si_*(\uh^{-1})=\sum_{n=0}^{\infty} \si(u^{-1})\otimes(1-\si(u)\otimes\si(u^{-1}))^{\otimes n}\ \in M_\infty(\Tt\Ec)^+\ .
$$
Under the hypotheses of Theorem \ref{tres}, the evaluation of $\Ind_E([u])$ on a cyclic cohomology class $[\varphi]\in HP_\top^{0}(\cinfc(B)\cp G)_{[O]}$ localized at $O$ is computed by the residue formula
\be
\langle [\tau_\varphi] , \Ind_E ([u]) \rangle =   \varphi\circ\Res\#\tr\big( \si_*(\uh^{-1}) [\ln Q, \si_*(\uh)] \big) \label{ctr} 
\ee
where $Q\in \cinf(B,\CL^1(M))$ is any elliptic operator as above.

\section{Geometric cocycles}\label{sgeo}

For any Lie groupoid $G\rightrightarrows B$ we shall construct in a geometric way cyclic cohomology classes $[\varphi]\in HP_\top^\bullet(\cinfc(B)\cp G)_{[O]}$ localized at $\Ad$-invariant isotropic submanifolds $O\subset G$. Our approach is intended to combine all the already known cohomologies giving rise to cyclic cocycles localized at the unit submanifold $O=B$: the cohomology $H^\bullet(BG)$ of the classifying space of an \'etale groupoid \cite{C94}, the differentiable cohomology $H^\bullet_d(G)$ of any Lie groupoid \cite{PPT}, the Gelfand-Fuchs cohomology appearing in the transverse geometry of a foliation \cite{C83}, etc..., and generalize them in order to include cyclic cohomology classes localized at any isotropic submanifold $O$. We start with a definition extending the notion of Cartan connection on a Lie groupoid; note that all the formalism developed here could certainly be reinterpreted into the language of multiplicative forms of \cite{CSS}. 

\begin{definition}
A \emph{connection} on a Lie groupoid $(r,s):\Gamma\rightrightarrows E$ is a Lie subgroupoid $\Phi \rightrightarrows F$ of the tangent groupoid $(r_*,s_*):T\Gamma\rightrightarrows TE$ such that
\begin{itemize}
\item $\Phi$ a vector subbundle of $T\Gamma$ simultaneously transverse to $\ker r_*$ and $\ker s_*$, in the sense that $\Phi\cap \ker r_*=\Phi\cap \ker s_*$ is the zero section of $T\Gamma$,
\item $F$ is a vector subbundle of $TE$, of the same rank as $\Phi$.
\end{itemize} 
The connection is \emph{flat} if $(\Phi,F)$ are integrable as subbundles of the tangent bundles $(T\Gamma,TE)$.
\end{definition}
By the transversality hypothesis, the tangent maps $r_*$ and $s_*$ send the fibers of $\Phi$ isomorphically onto the fibers of $F$. If the connection is flat, $(\Phi,F)$ define regular foliations on the manifolds $(\Gamma,E)$ respectively. The leaves of $\Gamma$ are locally diffeomorphic to the leaves of $E$ under both maps $r$ and $s$.\\
When $\Phi$ has maximal rank equal to the dimension of $E$, then $\Phi\oplus \ker r_*=\Phi\oplus \ker s_*=T\Gamma$ and $F=TE$. In that case, a flat connection defines a foliation of the manifold $\Gamma$ which is simultaneously transverse to the submersions $r$ and $s$ and of maximal dimension. This foliation in turn determines a set of local bisections of $\Gamma$. Since $\Phi$ is a subgroupoid of $T\Gamma$, one can check that these local bisections form a Lie pseudogroup with respect to the composition product on $\Gamma$. This pseudogroup acts on the manifold $E$ by local diffeomorphisms: any small enough local bisection provides a diffeomorphism from its source set in $E$ to its range set in $E$. Hence, any morphism $\gamma\in\Gamma$ can be extended to a local diffeomorphism from an open neighborhood of $s(\gamma)$ to an open neighborhood of $r(\gamma)$, in a way compatible with the composition of morphisms in $\Gamma$. \\
When $\Phi$ has rank $<\dim E$, the corresponding local bisections are only defined above each leaf of the foliation $F$ on $E$. This means that any morphism $\gamma\in\Gamma$ can be extended to a local diffeomorphism from a small open subset of the leaf containing $s(\gamma)$ to a small open subset of the leaf containing $r(\gamma)$. In general there are topological or geometric obstructions to the existence of a connection (flat or not) on a Lie groupoid. Here are some basic examples:
\begin{example}\textup{An \'etale groupoid $\Gamma\rightrightarrows B$ has a unique connection of maximal rank $T\Gamma$ which is always flat. The corresponding pseudogroup of local bisections is thus the pseudogroup of all local bisections of $\Gamma$. Hence, any morphism $\gamma\in \Gamma$ in an \'etale groupoid determines a local diffeomorphism from a neighborhood of the source $s(\gamma)$ to a neighborhood of the range $r(\gamma)$ in a unique way.}
\end{example}
\begin{example}\textup{Let $G$ be a Lie group acting on a manifold $B$ by \emph{global} diffeomorphisms. Then the action groupoid $\Gamma=B\rtimes G$ is endowed with the canonical flat connection $\ker(\mathrm{pr}_*)\subset T\Gamma$, where $\mathrm{pr}:\Gamma\to G$ is the projection. In this case the pseudogroup of local bisections determined by the connection is precisely the group $G$.}
\end{example}

A flat connection $\Phi\rightrightarrows F$ on a Lie groupoid $\Gamma\rightrightarrows E$ leads to the $\Gamma$-equivariant \emph{leafwise} cohomology of $E$. Indeed at any point $\gamma\in\Gamma$, the fiber $\Phi_\gamma$ is canonically isomorphic, as a vector space, to the fibers $s_*(\Phi_\gamma)=F_{s(\gamma)}$ and $r_*(\Phi_\gamma)=F_{r(\gamma)}$. Hence any $\gamma$ yields a linear isomorphism from the fiber $F_{s(\gamma)}$ to the fiber $F_{r(\gamma)}$, in a way compatible with the composition law in $\Gamma$. This means that the vector bundle $F$ over $E$ is a $\Gamma$-bundle. Here the flatness of the connection is not used. In the same way, the dual bundle $F^*$ is also a $\Gamma$-bundle. In general the cotangent bundle $T^*E$ may not carry any action of $\Gamma$, but the line bundle $\La^{\max} A\Gamma\otimes\La^{\max} T^*E$ always does. As usual $A\Gamma$ is the Lie algebroid of $\Gamma$. If $F_\perp = TE/F$ denotes the normal bundle to $F$, then one has a canonical isomorphism $\La^{\max} T^*E \cong \La^{\max} F_\perp^*\otimes\La^\top F^*$. Since $\Gamma$ acts on the line bundle $\La^{\max} F$ through the connection $\Phi$, one sees that the tensor product $\La^{\max} A\Gamma \otimes \La^{\max} T^*E \otimes \La^{\max} F \cong \La^{\max} A\Gamma \otimes \La^{\max} F_\perp^*$ is a $\Gamma$-bundle.  We define the bicomplex of $\Gamma$-equivariant leafwise differential forms $C^\bullet(\Gamma,\La^\bullet F) = (C^n(\Gamma,\La^m F))_{n\geq 0,m\geq 0}$ where 
\be
C^n(\Gamma,\La^m F) = \cinf(\Gamma^{(n)}, s^*(|\La^{\max} A\Gamma|\otimes |\La^{\max} F_\perp^*|\otimes \La^mF^*))
\ee
is the space of smooth sections of the vector bundle $|\La^{\max} A\Gamma|\otimes |\La^{\max} F_\perp^*|\otimes \La^mF^*$,  pulled back on the manifold $\Gamma^{(n)}$ of composable $n$-tuples by the source map $s:\Gamma^{(n)}\to E$, $(\gamma_1,\ldots,\gamma_n)\mapsto s(\gamma_n)$. In particular for $n=0$, $C^0(\Gamma,\La^mF) = \cinf(E, |\La^{\max} A\Gamma|\otimes |\La^{\max} F_\perp^*|\otimes \La^mF^*)$ is the space of leafwise $m$-forms on $E$ twisted by the line bundle $|\La^{\max} A\Gamma|\otimes |\La^{\max} F_\perp^*|$. Note that $|\La^{\max} F_\perp^*|$ is the bundle of $1$-densities transverse to the foliation $F$ on $E$. The first differential $d_1: C^n(\Gamma, \La^m F) \to C^{n+1}(\Gamma,\La^m F)$ on this bicomplex is the usual differential computing the groupoid cohomology with coefficients in a $\Gamma$-bundle. It is given on any cochain $c\in C^n(\Gamma,\La^mF)$ by
\beq
(d_1c) (\gamma_1,\ldots,\gamma_{n+1}) &=& c (\gamma_2,\ldots,\gamma_{n+1}) + \sum_{i=1}^n (-1)^i c (\gamma_1,\ldots,\gamma_i\gamma_{i+1},\ldots,\gamma_{n+1}) \non\\
&& + (-1)^{n+1}\, U_{\gamma_{n+1}}^{-1}\cdot(c (\gamma_1,\ldots,\gamma_n)) \non
\eeq
for all $(\gamma_1,\ldots,\gamma_{n+1})\in \Gamma^{(n+1)}$. The last term of the r.h.s. denotes the action of the linear isomorphism $U_{\gamma_{n+1}}^{-1}:(|\La^{\max} A\Gamma|\otimes |\La^{\max} F_\perp^*|\otimes \La^mF^*)_{s(\gamma_n)} \to (|\La^{\max} A\Gamma|\otimes |\La^{\max} F_\perp^*|\otimes \La^mF^*)_{s(\gamma_{n+1})}$ on $c (\gamma_1,\ldots,\gamma_n)$. The second differential $d_2 : C^n(\Gamma,\La^mF)\to C^n(\Gamma,\La^{m+1}F)$ comes from the leafwise de Rham differential $d_F: \cinf(E,\La^mF)\to \cinf(E,\La^{m+1}F)$ on the foliated manifold $E$. Indeed we first observe that the foliation $F$ on $E$ defines the sheaf of holonomy-invariant sections of the line bundle of transverse $1$-densities $|\La^{\max} F_\perp^*|$, which in turn induces a canonical foliated connexion $d_{F_\perp}: \cinf(E,|\La^{\max} F_\perp^*|)\to \cinf(E,|\La^{\max} F_\perp^*|\otimes F^*)$. This connexion is flat in the usual sense $(d_{F_\perp})^2=0$. In the same way, viewing $r^*A\Gamma$ as a subbundle of $T\Gamma$, the foliation $\Phi$ on $\Gamma$ defines the sheaf of holonomy-invariant sections of the line bundle $r^*|\La^{\max} A\Gamma|$ which descends to the line bundle $|\La^{\max} A\Gamma|$ over $E$, and subsequently defines a flat foliated connection $d_{A\Gamma}:\cinf(E,|\La^{\max} A\Gamma|)\to \cinf(E,|\La^{\max} A\Gamma|\otimes F^*)$. The sum $d_{A\Gamma} + d_{F_\perp} + d_F$ is a $\Gamma$-equivariant operator on the space of sections of the vector bundle $|\La^{\max} A\Gamma|\otimes |\La^{\max} F_\perp^*|\otimes \La^\bullet F^*$ which squares to zero. Finally we extend the foliation $\Phi$ on $\Gamma$ to a foliation $\Phi^{(n)}$ on $\Gamma^{(n)}$ in such a way that all projection maps $p_i:\Gamma^{(n)} \to \Gamma$, $(\gamma_1,\ldots,\gamma_n)\mapsto \gamma_i$ induce vector space isomorphisms $(p_i)_*:\Phi^{(n)}_{(\gamma_1,\ldots,\gamma_n)}\to \Phi_{\gamma_i}$. Since $\Phi$ is a subgroupoid of $\Gamma$ the foliation $\Phi^{(n)}$ exists and is unique. Moreover the source map $s:\Gamma^{(n)} \to E$ is a local diffeomorphism from the leaves of $\Phi^{(n)}$ to the leaves of $F$. We use this local identification to lift $d_{A\Gamma} + d_{F_\perp} + d_F$ to an operator $s^*(d_{A\Gamma} + d_{F_\perp} + d_F)$ on the space of sections of the vector bundle $s^*(|\La^{\max} A\Gamma|\otimes |\La^{\max} F_\perp^*|\otimes \La^\bullet F^*)$ over $\Gamma^{(n)}$. For any cochain $c\in C^n(\Gamma,\La^mF)$ we set
$$
(d_2c) = (-1)^n s^*(d_{A\Gamma} + d_{F_\perp} + d_F)c \ .
$$
Since $d_{A\Gamma} + d_{F_\perp} + d_F$ is $\Gamma$-equivariant, the two differentials on $C^\bullet(\Gamma,\La^\bullet F^*)$ anticommute: $d_1d_2+d_2d_1=0$. The $\Gamma$-equivariant leafwise cohomology of $E$ is by definition the cohomology of the total complex obtained from this bicomplex. If $\eta:E\to N$ is a submersion, we define $C^\bullet_\eta(\Gamma,\La^\bullet F^*)$ as the subcomplex of cochains having \emph{proper} support with respect to $\eta$.

\begin{definition}
Let $G\rightrightarrows B$ be a Lie groupoid and let $O$ be an $\Ad$-invariant isotropic submanifold of $G$. A \emph{geometric cocycle localized at} $O$ is a quadruple $(N,E,\Phi,c)$ where
\begin{itemize}
\item $N\stackrel{\nu}{\longrightarrow} B$ is a surjective submersion. Hence the pullback groupoid $\nu^* G\rightrightarrows N$ acts on the isotropic submanifold $\nu^*O\subset\nu^*G$ by the adjoint action.
\item $E\stackrel{\eta}{\longrightarrow}\nu^*O$ is a $\nu^*G$-equivariant submersion. Any element $\gamma\in\nu^*O$, viewed in $\nu^*G$, is required to act by the identity on its own fiber $E_\gamma$.
\item $\Phi\rightrightarrows F$ is a flat connection on $\Gamma=E\cp \nu^*G$, with $F$ oriented. The canonical section $E\to\nu^*O\to\Gamma$, restricted to a leaf of $F$, is a leaf of $\Phi$.
\item $c\in C_\eta^{\bullet}(\Gamma, \La^\bullet F^*)$ is a total cocycle with proper support relative to the submersion $E\stackrel{\eta}{\longrightarrow}\nu^*O$, assumed \emph{normalized} in the sense that
$$
c(\gamma_1,\ldots,\gamma_n) = 0 \quad \mbox{whenever} \quad \gamma_1\ldots\gamma_n = \eta(s(\gamma_n))\ .
$$ 
\end{itemize}
$(N,E,\Phi,c)$ is called \emph{proper} if $E$ is a proper $\nu^*G$-manifold.
\end{definition}

From now on let $G\rightrightarrows B$ be a fixed Lie groupoid. We shall associate to any geometric cocycle $(N,E,\Phi,c)$ localized at an isotropic submanifold $O\subset G$, a periodic cyclic cohomology class of $\cinfc(G)$ localized at $O$. This requires a number of steps. On the action groupoid $\Gamma = E\cp \nu^*G$ we first define the convolution algebra
\be 
\Gc = \cinf_p(E, \La^\bullet F^* )\cp \Gamma
\ee
which is the space $\cinf_p(\Gamma, r^*(\La^\bullet F^* \otimes |\La^{\max} A^*\Gamma|))$ of properly supported sections of the vector bundle $|\La^{\max} A^*\Gamma|\otimes \La^\bullet F^*$ pulled back by the rank map $r:\Gamma\to E$. By definition the groupoid $\Gamma$ acts on the vector bundle $\La^\bullet F^*$, but there is no such action on the density bundle $|\La^{\max} A^*\Gamma|$. For any $\gamma\in \Gamma$ we thus define the linear isomorphism $U_\gamma: \La^\bullet F^*_{s(\gamma)} \to \La^\bullet F^*_{r(\gamma)}$ leaving the space $|\La^{\max} A^*\Gamma|_{s(\gamma)}$ untouched. The convolution product of two elements $\al_1,\al_2\in \Gc$ is then given by 
\be 
(\al_1\al_2)(\gamma) = \int_{\gamma_1\gamma_2=\gamma} \al_1(\gamma_1) \wedge U_{\gamma_1} \al_2(\gamma_2)\ ,\qquad \forall \gamma\in\Gamma \ , \label{convg}
\ee
where $\al_1(\gamma_1) \wedge U_{\gamma_1} \al_2(\gamma_2) \in \La^\bullet F^*_{r(\gamma_1)} \otimes |\La^{\max} A^*\Gamma|_{r(\gamma_1)} \otimes |\La^{\max} A^*\Gamma|_{r(\gamma_2)}$ involves the exterior product of leafwise differential forms, and the integral is taken against the $1$-density $|\La^{\max} A^*\Gamma|_{r(\gamma_2)}$ while $r(\gamma_1)=r(\gamma)$ remains fixed. Hence $(\al_1\al_2)(\gamma)$ really defines an element of the fiber $\La^\bullet F^*_{r(\gamma)} \otimes |\La^{\max} A^*\Gamma|_{r(\gamma)}$. Note that the subalgebra of $\Gc$ consisting only in the zero-degree foliated differential forms coincides with the usual convolution algebra $\cinf_p(E)\cp\Gamma \cong \cinf_p(\Gamma,r^*|\La^{\max} A^*\Gamma|)$ of (properly supported) scalar functions over the groupoid $\Gamma$. Now let $d_{A^*\Gamma}:\cinf(E,|\La^{\max} A^*\Gamma|)\to \cinf(E,F^*\otimes |\La^{\max} A^*\Gamma|)$ be the flat foliated connection dual to $d_{A\Gamma}$, induced as before by the connection $\Phi$. Combining $d_{A^*\Gamma}$ with the leafwise de Rham differential $d_F:\cinf(E,\La^m F^*)\to \cinf(E,\La^{m+1} F^*)$, we obtain a total differential $d_\Phi$ on $\Gc$,
$$
d_\Phi\al = r^*(d_F + d_{A^*\Gamma}) \al \ ,\qquad \al\in \Gc
$$
where $d_F + d_{A^*\Gamma}$, acting leafwise on $E$, is lifted to a leafwise operator on $\Gamma$ through the local identification between the leaves of $F$ and the leaves of $\Phi$ provided by the rank map. Then $d_\Phi$ satisfies the graded Leibniz rule $d_\Phi(\al_1\al_2) = (d_\Phi\al_1)\al_2 + (-1)^{|\al_1|}\al_1 (d_\Phi\al_2)$, where $|\al_1|$ is the degree of the differential form $\al_1$. Hence $(\Gc,d_\Phi)$ is a differential graded (DG) algebra.\\
The periodic cyclic cohomology of a DG algebra is defined in complete analogy with the usual case, simply by adding the extra differential and taking care of the degrees of the elements in the algebra. Hence the space of noncommutative differential forms $\Om\Gc$ is the same as in the ungraded case, but the degree of an $n$-form $\om = \al_0d\al_1\ldots d\al_n$ is now $|\om| = n+ |\al_0| +\ldots +|\al_n|$. The differential $d_\Phi$ is uniquely extended to a differential on the graded algebra $\Om\Gc$, in such a way that it anticomutes with $d$. Hence we have
$$
d_\Phi(\al_0d\al_1\ldots d\al_n) = (d_\Phi\al_0)d\al_1\ldots d\al_n - (-1)^{|\al_1|} \al_0d(d_\Phi\al_1)\ldots d\al_n + \ldots
$$
The Hochschild boundary is as usual $b(\om\al) = (-1)^{|\om|}[\om,\al]$ for any $\om\in\Om\Gc$ and $\al\in\Gc$, where the commutator is the \emph{graded} one. Then $d_\Phi$ anticommutes with $b$, with the Karoubi operator $\kappa = 1-(bd+db)$, and with Connes'operator $B=(1+\kappa+\ldots+\kappa^n)d$ on $\Om^n\Gc$. The periodic cyclic cohomology of $\Gc$ is therefore defined as the cohomology of the complex $\hom(\Omh\Gc,\cc)$ with boundary map the transposed of $b+B+d_\Phi$.\\
We now define a linear map $\la:C^\bullet_\eta(\Gamma,\La^\bullet F^*)\to \hom(\Omh\Gc,\cc)$ which, once restricted to \emph{normalized} cochains, will behave like a chain map. For any cochain $c\in C^n_\eta(\Gamma,\La^m F^*)$ set
\beq
\lefteqn{\la(c) (\al_0d\al_1\ldots d\al_n) = }\non\\
&&\int_{(\gamma_1,\ldots,\gamma_n)\in\Gamma^{(n)}} \al_0(\gamma_0) \wedge U_{\gamma_0}\al_1(\gamma_1) \ldots \wedge U_{\gamma_0\ldots\gamma_{n-1}}\al_n(\gamma_n) \wedge c(\gamma_1,\ldots,\gamma_n) \non
\eeq
where $\gamma_0$ is defined as a function of the $n$-tuple $(\gamma_1\ldots\gamma_n)$ by the localization condition $\gamma_0\ldots\gamma_n = \eta(s(\gamma_n))$. We use the vector space isomorphism $|\La^{\max} A^*\Gamma|_{r(\gamma_0)}\otimes |\La^{\max} A\Gamma|_{s(\gamma_n)} \cong \cc$ to view the wedge product under the integral as an element of the fiber $|\La^{\max} A^*\Gamma|_{r(\gamma_1)}\otimes \ldots \otimes |\La^{\max} A^*\Gamma|_{r(\gamma_n)} \otimes \La^{m+|\al|} F^*_{s(\gamma_n)}$, with $|\al|=|\al_0|+\ldots +|\al_n|$. Since $F$ is oriented,  the integrand defines a $1$-density which can be integrated over the manifold $\Gamma^{(n)}$ when the leafwise degree $m+|\al|$ matches the rank of $F$; otherwise the integral is set to zero.

\begin{lemma}
Let $c\in C^\bullet_\eta(\Gamma,\La^\bullet F^*)$ be a \emph{normalized} cochain in the sense that $c(\gamma_1,\ldots,\gamma_n) = 0$ whenever $\gamma_1\ldots\gamma_n = \eta(s(\gamma_n))$. Then the periodic cyclic cochain $\la(c)\in \hom(\Omh\Gc,\cc)$ verifies the identities 
$$
\la(c)\circ b = \la(d_1c)\ ,\qquad \la(c)\circ d  = 0\ ,\qquad  \la(c) \circ d_\Phi = \la(d_2c)\ .
$$
Hence if $c$ is a normalized total cocycle, $\la(c)$ is a $\kappa$-invariant periodic cyclic cocycle over the DG algebra $\Gc$.
\end{lemma}
{\it Proof:} The three identities are routine computations. Since $\la(c)\circ d=0$ for any normalized $c$, one has $\la(c)\circ B=0$. As a consequence $(d_1+d_2)c = 0$ implies $\la(c)\circ(b+B+d_\Phi)=0$, thus $\la(c)$ is a periodic cyclic cocycle of the DG algebra $\Gc$. Moreover for any normalized cocycle $c$,  
$$
\la(c)\circ(1-\kappa) = \la(c)\circ(db+bd) = \la(d_1c)\circ d = -\la(d_2c)\circ d = -\la(c)\circ d_\Phi d = 0
$$ 
since $d_\Phi$ and $d$ anticommute, which shows that the periodic cyclic cocycle $\la(c)$ is $\kappa$-invariant. \cqfd\\

The $\kappa$-invariance of the cocycle $\la(c)$ means that the latter can as well be interpreted as an $X$-complex cocycle for certain DG algebra extensions of $(\Gc,d_\Phi)$. The $X$-complex of any associative DG algebra $(\Hc,d)$ is defined in analogy with the usual case by
$$
X(\Hc,d) : \Hc \rightleftarrows \Om^1\Hc_\nat\ ,
$$
where $\Om^1\Hc_\nat$ is the quotient of $\Om^1\Hc$ by the subspace of \emph{graded} commutators $[\Hc,\Om^1\Hc]=b\Om^2\Hc$. Since $d$ acts on $\Om\Hc$ and anticommutes with the Hochschild operator $b$, it descends to a well-defined differential on $X(\Hc,d)$ and anticommutes with the usual $X$-complex boundary maps $\nat\dd:\Hc \to \Om^1\Hc_\nat$ and $\bb:\Om^1\Hc_\nat\to\Hc$. We always endow the $X$-complex of a DG algebra with the total boundary operator $(\nat\dd\oplus\bb) +d$. Now take $\Hc$ as the direct sum $\Hc = \bigoplus_{n\geq 1} \Hc_n$, where
\be  
\Hc_n = \cinf_p(\Gamma^{(n)}, r^*_1\La^\bullet F^* \otimes r^*_1|\La^{\max} A^*\Gamma|\otimes \ldots \otimes r^*_n|\La^{\max} A^*\Gamma|)
\ee 
and $r_i:\Gamma^{(n)}\to E$ is the rank map $(\gamma_1,\ldots\gamma_n)\mapsto r(\gamma_i)$. The component $\Hc_1$ is isomorphic, as a vector space, to $\Gc$. The product of two homogeneous elements $\al_1\in \Hc_{n_1}$ and $\al_2\in\Hc_{n_2}$ is the element $\al_1\al_2\in \Hc_{n_1+n_2}$ defined by
$$
(\al_1\al_2)(\gamma_1,\ldots,\gamma_{n_1+n_2}) = \al_1(\gamma_1,\ldots,\gamma_{n_1}) \wedge U_{\gamma_1\ldots\gamma_{n_1}}\al_2(\gamma_{n_1+1},\ldots,\gamma_{n_1+n_2})
$$
We equip $\Hc$ with a grading by saying that an element $\al \in \cinf_p(\Gamma^{(n)}, r^*_1\La^m F^* \otimes r^*_1|\La^{\max} A^*\Gamma|\otimes \ldots \otimes r^*_n|\La^{\max} A^*\Gamma|)$ has degree $|\al| = m$. A differential $d$ of degree $+1$ on $\Hc$ is then defined by combining the leafwise differential $d_F$ with the flat connections on the density bundles $\La^{\max}|A^*\Gamma|$:
$$
d \al = (r_1^*(d_F) + r_1^*(d_{A^*\Gamma}) + \ldots + r_n^*(d_{A^*\Gamma}))\al
$$
for any such $\al\in \Hc_n$. Obviously $d$ satisfies the graded Leibniz rule. By construction $(\Hc,d)$ is an extension of $(\Gc,d_\Phi)$. Indeed a surjective DG algebra morphism $m:\Hc\to\Gc$ is defined as follows: the image of any $\al\in\Hc_n$ is the element $m(\al)\in \Gc$ given by the $(n-1)$-fold integral
$$
(m(\al))(\gamma) = \int_{\gamma_1\ldots\gamma_n = \gamma} \al(\gamma_1,\ldots,\gamma_n)\ ,\qquad \forall\gamma\in\Gamma
$$
where the integral is taken against the density bundle $r^*_2|\La^{\max} A^*\Gamma|\otimes \ldots \otimes r^*_n|\La^{\max} A^*\Gamma|$. Let $\Ic$ be the kernel of the multilication map. In fact $\Hc$ is a kind of localized version of the tensor algebra extension $T\Gc$ of $\Gc$. The latter is a graded algebra, the degree of a tensor $\al_1\otimes\ldots \otimes\al_n$ being the sum of the degrees of the factors. $T\Gc$ is endowed with a differential $d_\Phi$ making the multiplication map $T\Gc\to\Gc$ a morphism of DG algebras:
$$
d_\Phi(\al_1\otimes\ldots \otimes\al_n) = (d_\Phi\al_1)\otimes\ldots \otimes\al_n + (-1)^{|\al_1|} \al_1\otimes(d_\Phi\al_2)\otimes\ldots \otimes\al_n + \ldots
$$
A straightforward adaptation of the Cuntz-Quillen theory to the DG case shows that the $X$-complex $X(\Th\Gc,d_\Phi)$ is quasi-isomorphic to the $(b+B+d_\Phi)$-complex of noncommutative differential forms $\Omh\Gc$, under the identification of pro-vector spaces $X(\Th\Gc)\cong\Omh\Gc$ taking the rescaling factors $(-1)^{[n/2]}[n/2]!$ into account. Since the periodic cyclic cocycle $\la(c)$ is $\kappa$-invariant, it can as well be viewed as an $X$-complex cocycle:
$$
\la'(c)\in\hom(X(\Th\Gc,d_\Phi),\cc)\ .
$$
In fact this cocycle descends to a cocycle over $X(\Hch,d)$, where $\Hch$ denotes the $\Ic$-adic completion of $\Hc$. Indeed we note that the canonical \emph{linear} inclusion $\Gc\hookrightarrow \Hc$, which identifies $\Gc$ and the vector subspace $\Hc_1$, commutes with the differential. The universal property of the tensor algebra then implies the existence of a DG algebra homomorphism $T\Gc\to\Hc$, mapping the DG ideal $J\Gc$ to $\Ic$, whence a chain map $X(\Th\Gc,d_\Phi)\to X(\widehat{\Hc},d)$. Observe that the homomorphism $T\Gc\to\Hc$ has dense range. The following lemma is obvious.

\begin{lemma}
Let $c\in C^\bullet_\eta(\Gamma,\La^\bullet F^*)$ be a normalized cocycle. Then $\la'(c)$, viewed as a cocycle in $\hom(X(\Th\Gc,d_\Phi),\cc)$, factors through a unique continuous cocycle $\la'(c) \in \hom(X(\Hch_\top,d),\cc)$. \cqfd
\end{lemma}

The last step is the construction of an homomorphism from the convolution algebra of compactly supported functions $\cinfc(B)\cp G$ to the convolution algebra of \emph{properly} supported functions $\cinf_p(E)\cp\Gamma$. Indeed, the pullback groupoid $\nu^*G\rightrightarrows N$ is Morita equivalent to $G\rightrightarrows B$, so using a cut-off function on $N$ we realize the equivalence by an homomorphism
$$
\cinfc(B)\cp G \stackrel{\sim}{\longrightarrow} \cinfc(N)\cp \nu^*G\ .
$$
Next, by definition the submersion $E\stackrel{\eta}{\to} \nu^*O \stackrel{s}{\to} N$ is $\nu^*G$-equivariant. Hence any $\gamma\in\Gamma=E\cp\nu^*G$ determines a unique $g\in \nu^*G$, and the resulting map $\Gamma\to \nu^*G$ is a smooth morphism of Lie groupoids. The latter identifies the preimage of the source map $\Gamma_{s(\gamma)}$ in $\Gamma$ with the preimage of the source map $\nu^*G_{s(g)}$ in $\nu^*G$, whence a canonical vector space isomorphism between the fibers $A\Gamma_{r(\gamma)}$ and $A(\nu^*G)_{r(g)}$ of the corresponding Lie algebroids. We conclude that any smooth compactly supported section of $|\La^{\max} A^*(\nu^*G)|$ over $\nu^*G$ can be canonically pulled back to a smooth properly supported section of $|\La^{\max} A^*\Gamma|$ over $\Gamma$, and this results in an homomorphism of convolution algebras $\cinfc(N)\cp\nu^*G \to\cinf_p(E)\cp\Gamma$. By composition with the Morita equivalence, we obtain the desired homomorphism
\be 
\rho: \cinfc(B) \cp G \to \cinf_p(E)\cp\Gamma\ .
\ee
The latter depends on the choice of homomorphism realizing the Morita equivalence. However this dependence will disappear in cohomology. Now put $\Cc = \cinfc(B)\cp G$. Since $\cinf_p(E)\cp\Gamma$ is the subalgebra of $\Gc$ consisting of zero-degree foliated differential forms, we regard $\rho$ as an algebra homomorphism $\Cc\to\Gc$. The latter lifts to a linear map $\iota:\Cc\to\Hc$ after composition by the canonical linear inclusion $\Gc=\Hc_1\hookrightarrow\Hc$. The diagram of extensions
\be
\vcenter{\xymatrix{0 \ar[r] & J\Cc \ar[r] \ar[d]_{\rho_*} & T\Cc \ar[r] \ar[d]_{\rho_*} & \Cc \ar[r] \ar@{.>}[dl]^{\iota} \ar[d]^{\rho} & 0 \\
0\ar[r] & \Ic \ar[r]  & \Hc \ar[r] & \Gc \ar[r]   & 0 }}
\ee
thus allows to extend $\rho_*$ to an homomorphism of pro-algebras $\Th\Cc\to\Hch$. 

\begin{lemma}\label{lchi}
The linear map $\chi(\rho_*,d) \in \hom(\Omh\Th\Cc , X(\Hch_\top,d))$ defined on any $n$-form $\hat{c}_0 \dd\hat{c}_1\ldots \dd\hat{c}_n$ by 
\beq 
\lefteqn{\chi(\rho_*,d)(\hat{c}_0 \dd\hat{c}_1\ldots \dd\hat{c}_n)=} \\
&& \frac{1}{(n+1)!} \sum_{i=0}^n (-1)^{i(n-i)} d\rho_*(\hat{c}_{i+1}) \ldots d\rho_*(\hat{c}_n) \, \rho_*(\hat{c}_0) \, d\rho_*(\hat{c}_1) \ldots d\rho_*(\hat{c}_i) \non \\
&& + \frac{1}{n!} \sum_{i=1}^n \nat\big(\rho_*(\hat{c}_0) \, d\rho_*(\hat{c}_1) \ldots \dd \rho_*(\hat{c}_i) \ldots d\rho_*(\hat{c}_n) \big) \non
\eeq
is a chain map from the $(b+B)$-complex of noncommutative differential forms to the DG $X$-complex. Moreover the cohomology class of $\chi(\rho_*,d)$ in the $\hom$-complex $\hom(\Omh\Th\Cc,X(\Hch_\top,d))$ is independent of any choice concerning the homomorphism $\rho$.
\end{lemma}
{\it Proof:} A routine computation shows that $\chi(\rho_*,d)$ is a chain map. The independence of its cohomology class upon the choice of homomorphism $\rho$ is a classical homotopy argument using $2\times 2$ rotation matrices. \cqfd\\

\begin{proposition}
Let $G\rightrightarrows B$ be a Lie groupoid and let $O\subset G$ be an $\ad$-invariant isotropic submanifold. Any geometric cocycle $(N,E,\Phi,c)$ localized at $O$ defines a class $[N,E,\Phi,c]\in HP_\top^{\bullet}(\cinfc(B)\cp G)_{[O]}$, represented by the composition of chain maps
$$
\xymatrix{X(\Th\Cc) \ar[r]^{\gamma} & \Omh\Th\Cc \ar[r]^{\chi(\rho_*,d)\quad} & X(\Hch_\top,d) \ar[r]^{\qquad \la'(c)} & \cc }
$$
where $\gamma$ is the generalized Goodwillie equivalence, $\Cc=\cinfc(B)\cp G$ is the convolution algebra of $G$, and $(\Hch_\top,d)$ is the DG pro-algebra constructed above from the geometric cocycle.  \cqfd
\end{proposition}

\section{Localization at units}\label{sloc}

Let $G\rightrightarrows B$ be a Lie groupoid and $(N,E,\Phi,c)$ a geometric cocycle localized at units. Hence $\nu:N\to B$ is a surjective submersion, $\eta:E\to N$ is a $\nu^*G$-equivariant submersion, $\Phi\rightrightarrows F$ is a flat connection on the action groupoid $\Gamma = E\cp \nu^*G$, and $c\in C_\eta^\bullet(\Gamma,\La^\bullet F^*)$ is a normalized cocycle. \emph{Throughout this section we assume that $E$ is a proper $\nu^*G$-manifold}. Let $\pi:M\to B$ be a $G$-equivariant submersion. We make no properness hypothesis about the action of $G$ on $M$. Then the algebra bundle of vertical symbols $\CS_c(M)$ over $B$ is a $G$-bundle. Its pullback $\CS_c(N\times_BM)$ under the submersion $\nu$ is a bundle over $N$, whose fibers are isomorphic to the same algebras of vertical symbols. The pullback groupoid $\nu^*G \rightrightarrows N$ acts naturally on $\CS_c(N\times_BM)$. By hypothesis the action of $\nu^*G$ on $N$ also lifts to $E$, hence the pullback of $\CS_c(N\times_BM)$ under the submersion $\eta$ yields a $\Gamma$-bundle $\CS_c(E\times_BM)$ over $E$. The vector bundle $F\subset TE$ being also a $\Gamma$-bundle, we can form the convolution algebra
\be
\Oc = \cinf_p(E,\La^\bullet F^*\otimes\CS_c(E\times_B M))\cp\Gamma \ ,
\ee
which is a symbol-valued generalization of the algebra $\Gc$ of section \ref{sgeo}. The product on $\Oc$ is formally identical to (\ref{convg}), involving the algebra structure of the bundle $\La^\bullet F^*\otimes\CS_c(E\times_BM)$ together with the linear isomorphism $U_\gamma : (\La^\bullet F^*\otimes\CS_c(E\times_BM))_{s(\gamma)} \to (\La^\bullet F^*\otimes\CS_c(E\times_BM))_{r(\gamma)}$ for all $\gamma \in \Gamma$. The agebra $\Oc$ is naturally graded by the form degree in $\La^\bullet F^*$.\\ 
Let $(\id,\pi)_*$ be the tangent map of the submersion $(\id,\pi):E\times_BM\to E$. The preimage of the integrable subbundle $F\subset TE$ is an integrable subbundle $(\id,\pi)_*^{-1}(F)$ of $T(E\times_BM)$ defining a foliation on $E\times_BM$. Choose an horizontal distribution $H$ in this subbundle, that is, a decomposition $(\id,\pi)_*^{-1}(F) = H\oplus \ker(\id,\pi)_*$. By construction the groupoid $\Gamma$ acts on $(\id,\pi)_*^{-1}(F)$, and by properness we can even assume that $H$ is $\Gamma$-invariant if necessary. Let $\cinfc(E\times_BM)\to E$ be the bundle over $E$ whose fibers are smooth vertical functions with compact support. We can identify this bundle with $\PS^0_c(E\times_BM)\to E$, the polynomial vertical symbols of order 0. Combining the distribution $H$ with the leafwise de Rham differential $d_F:\cinf(E,\La^mF^*) \to \cinf(E,\La^{m+1}F^*)$ yields a ``foliated" connection on this bundle, in the sense of a linear map
$$
d_H\ :\ \cinf(E,\La^m F^*\otimes\cinfc(E\times_BM)) \to \cinf(E,\La^{m+1} F^*\otimes\cinfc(E\times_BM))
$$
which is a derivation of $\cinf(E)$-modules. In general the subbundle $H$ is not integrable and $d_H$ does not square to zero. Its curvature
$$
(d_H)^2 = \te \ \in \cinf(E,\La^2F^*\otimes \CS^1(E\times_BM))
$$
is a $\Gamma$-invariant leafwise 2-form over $E$ with values in vertical vector fields. By definition the bundle of vertical symbols $\CS_c(E\times_BM)$ acts by endomorphisms on the bundle $\cinfc(E\times_BM)$. Hence the graded commutator $\tilde{d}_H = [d_H,\ ]$ is a graded derivation on the algebra of sections $\cinf(E,\La^\bullet F^*\otimes\CS_c(E\times_BM))$, with curvature $(\tilde{d}_H)^2 = [\te,\ ]$. Combining further $\tilde{d}_H$ with the flat connection
$$
d_{A^*\Gamma}: \cinf(E,|\La^{\max}A^*\Gamma|) \to \cinf(E, F^*\otimes |\La^{\max}A^*\Gamma|)
$$ 
as in section \ref{sgeo}, we get a derivation (still denoted by $\tilde{d}_H$) on the algebra $\Oc$. Then $(\tilde{d}_H)^2$ still acts by the commutator $[\te,\ ]$, where $\te$ is viewed as a multiplier of $\Oc$. In order to deal with cyclic cohomology we construct an extension of this algebra. Define the vector space $\Pc = \bigoplus_{n\geq 1} \Pc_n$, where
$$
\Pc_n = \cinf_p(\Gamma^{(n)}, r^*_1(\La^\bullet F^*\otimes\CS_c(E\times_B M)) \otimes r^*_1|\La^{\max} A^*\Gamma|\otimes \ldots \otimes r^*_n|\La^{\max} A^*\Gamma|)
$$
and $r_i:\Gamma^{(n)}\to E$ is the rank map $(\gamma_1,\ldots\gamma_n)\mapsto r(\gamma_i)$. The component $\Pc_1$ is isomorphic, as a vector space, to $\Oc$. The product of two elements $\al_1\in \Pc_{n_1}$ and $\al_2\in\Pc_{n_2}$ is the element $\al_1\al_2\in \Pc_{n_1+n_2}$ defined by
$$
(\al_1\al_2)(\gamma_1,\ldots,\gamma_{n_1+n_2}) = \al_1(\gamma_1,\ldots,\gamma_{n_1}) \wedge U_{\gamma_1\ldots\gamma_{n_1}}\al_2(\gamma_{n_1+1},\ldots,\gamma_{n_1+n_2})
$$
We equip $\Pc$ with the grading induced by $\La^\bullet F$. Then $\Pc$ is a symbol-valued generalization of the algebra $\Hc$ of section \ref{sgeo}. One has a multiplication homomorphism $m:\Pc\to\Oc$, and the derivation $\tilde{d}_H$ on $\Oc$ extends in a unique way to a derivation on the algebra $\Pc$. We let $\Qc$ be the ideal $\ker(m)$ and denote as usual by $\Pch$ the $\Qc$-adic completion of $\Pc$. Hence $\Pch$ is a graded pro-algebra, endowed with a derivation $\tilde{d}_H$ of degree $1$.\\
Let $\Ac=\cinfc(S^*_\pi M)\cp G$. We want to construct an homomorphism from $\Th\Ac$ to $\Pch$. Using a cut-off function on the submersion $\nu:N\to B$, we know that the Morita equivalence between the groupoids $G$ and $\nu^*G$ is realized by an homomorphism $\cinfc(B)\cp G\to \cinfc(N)\cp\nu^* G$ of the corresponding convolution algebras. Using the same cut-off function, we get an homomorphism
$$
\Ec = \cinfc(B,\CL^0_c(M))\cp G \to \cinfc(N,\CL^0_c(N\times_BM))\cp \nu^*G\ .
$$
Then as done in section \ref{sgeo} the $\Gamma$-equivariant map $\eta:E\to N$ induces a pullback homomorphim $\cinfc(N,\CL^0_c(N\times_BM))\cp \nu^*G \to \cinf_p(E,\CL^0_c(E\times_BM))\cp \Gamma$. Taking further the projection of classical pseudodifferential operators onto formal symbols one is left with an homomorphism
$$
\mu\ :\ \Ec \to \cinf_p(E,\CS^0_c(E\times_BM))\cp \Gamma \subset \Oc\ ,
$$
sending the ideal $\Bc=\cinfc(B,\CL^{-1}_c(M))\cp G$ to the convolution algebra of symbols of negative order $\cinf_p(E,\CS^{-1}_c(E\times_BM))\cp \Gamma$. By the linear inclusion $\Oc\hookrightarrow \Pc_1$ it extends to an homomorphism $\mu_*:T\Ec\to \Pc$, sending the ideal $J\Ec$ to $\Qc$. Moreover the image of the ideal $\Rc = T(\Bc:\Ec)\subset T\Ec$ contains only symbols of order $\leq -1$ in $\Pc$. Therefore $\mu_*$ first extends to an homomorphism from $\varprojlim_n T\Ec/\Rc^n$ to $\Pc$, and then from $\Tt\Ec$ to $\Pch$. Composing with the homomorphism $\si_*:\Th\Ac\to\Tt\Ec$ of section \ref{sres}, one thus gets a new homomorphism
\be 
\si_*' = \mu_*\circ\si_* \ :\ \Th\Ac\to \Tt\Ec \to \Pch 
\ee
Now we twist the algebra $\Pch$ by adding an odd parameter $\epsilon$, with the property $\epsilon^2 = 0$. Let $\Pch[\epsilon]$ be the resulting $\zz_2$-graded algebra: it is linearly spanned by elements of the form $\al_0+\epsilon\al_1$ for $\al_0,\al_1\in \Pch$, with obvious multiplication rules. Choose a section $Q\in\cinf(E,\CL^1(E\times_BM))$ of vertical elliptic operators of order one over $E$, with symbol $q\in \cinf(E,\CS^1(E\times_BM))$. By properness we can even assume that $Q$ and $q$ are $\Gamma$-invariant if necessary. The superconnection
\be
\nabla = d_H + \epsilon \ln q\ ,
\ee
acting by graded commutators, is an odd derivation on $\Pch[\epsilon]$. Indeed $[\nabla,\al] = \tilde{d}_H\al +\epsilon[\ln q,\al]\in\Pch[\eps]$ for all $\al\in\Pch[\eps]$. Moreover, any derivative of the logarithmic symbol $\ln q$ being a classical symbol, the curvature of the superconnection 
$$
\nabla^2 = \te -\epsilon \tilde{d}_H\ln q \ \in \cinf(E,\La^\bullet F^*\otimes\CS(E\times_BM)[\epsilon])
$$
is a multiplier of $\Pch[\eps]$. Let $\Hch$ be the algebra constructed in section \ref{sgeo}. From the homomorphism $\si_*'$ and the superconnection $\nabla$ we construct a chain map $\chi^\Res(\si_*',\nabla) \in \hom(\Omh\Th\Ac , X(\Hch_\top[\epsilon], d)_{[E]})$ by means of a JLO-type formula \cite{JLO}. Since the complex $X(\Hch_\top[\epsilon],d)_{[E]}$ is localized at units, the Wodzicki residue of vertical symbols gives a linear map
\be 
\Res\ :\ X(\Pch[\epsilon]) \to X(\Hch_\top[\epsilon],d)_{[E]}\ . \label{res}
\ee
For any $n$-form $\ah_0 \dd\ah_1 \ldots\dd\ah_n \in\Om^n\Th\Ac$ we set
\beq
\lefteqn{\chi^\Res(\si_*',\nabla)(\ah_0\dd\ah_1\ldots\dd\ah_n)=}\label{JLO1}\\
&&\hspace{-1cm}\sum_{i=0}^n (-)^{i(n-i)} \int_{\Delta_{n+1}} \Res\big(e^{-t_{i+1}\nabla^2} [\nabla,\si'_{i+1}] \ldots  e^{-t_{n+1}\nabla^2} \si'_0 \,e^{-t_0\nabla^2} [\nabla,\si'_1] \ldots e^{-t_i\nabla^2}\big)dt \non\\
&+&  \sum_{i=1}^n \int_{\Delta_n} \Res\big(\nat \si'_0\, e^{-t_0\nabla^2} [\nabla,\si'_1] \ldots  e^{-t_{i-1}\nabla^2} \dd \si'_i \,e^{-t_i\nabla^2} \ldots [\nabla,\si'_n] e^{-t_n\nabla^2}\big)dt \non
\eeq
where $\si'_i = \si'_*(\ah_i)$ for all $i$, and $\Delta_n =\{(t_0,\ldots,t_n)\in [0,1]^n \ |\ t_0+\ldots+t_n = 1\}$ is the standard $n$-simplex. This formula makes sense because $\nabla^2 = \te -\epsilon \tilde{d}_H\ln q$ is nilpotent as a leafwise differential form of degree $\geq 1$ over $E$. Hence the products under the residue are well-defined elements of $X(\Pch[\epsilon])$, depending polynomially on the simplex variable $t$. The nilpotency of the curvature also implies that $\chi^\Res(\si'_*,\nabla)$ vanishes on $\Om^n\Th\Ac$ whenever $n>\dim F + 2$. Basic computations show that (\ref{JLO1}) are the components of a chain map from the $(b+B)$-complex $\Omh\Th\Ac$ to the $X$-complex of the DG algebra $\Hc_\top[\epsilon]$ localized at $E$. Now $\chi^\Res(\si'_*,\nabla)$ may be expanded as a sum of terms which do not contain $\epsilon$, plus terms exactly proportional to $\epsilon$. We define the cocycle $\chi^\Res(\si'_*,d_H,\ln q) \in \hom(\Omh\Th\Ac , X(\Hch_\top, d)_{[E]})$ as the coefficient of $\epsilon$ in the latter expansion, or equivalently as the formal derivative 
\be 
\chi^\Res(\si'_*,d_H,\ln q) = \frac{\d}{\d\epsilon}\, \chi^\Res(\si'_*,d_H+\epsilon\ln q)\ .
\ee
By classical Chern-Weil theory, higher transgression formulas show that the cohomology class of the cocycle $\chi^\Res(\si'_*,d_H,\ln q)$ does not depend on the choice of connection $d_H$ and elliptic symbol $q$, and is a homotopy invariant of the homomorphism $\si_*'$.

\begin{proposition}\label{pcompo}
Let $\Ac$ be the convolution algebra of the action groupoid $S^*_\pi M\cp G$. Then image of the cyclic cohomology class of a proper geometric cocycle localized at units $(N,E,\Phi,c)$ under the excision map
$$
\pi^!_G : HP_\top^{\bullet}(\cinfc(B)\cp G)_{[B]} \to HP_\top^{\bullet+1}(\cinfc(S^*_\pi M)\cp G)_{[S^*_\pi M]}
$$ 
is the cyclic cohomology class over $\Ac$ represented by the chain map
\be 
\la'(c)\circ\chi^\Res(\si'_*,d_H,\ln q)\circ \gamma\ :\ X(\Th\Ac) \to \Omh\Th\Ac \to X(\Hch_\top,d)_{[E]} \to \cc \label{path}
\ee 
\end{proposition}
{\it Proof:} Since $(d_H)^2=\te\neq 0$, $(\Pch,\tilde{d}_H)$ is not a differential pro-algebra. We use a trick of Connes (\cite{C94} p.229) and add a multiplier $v$ of $\Pch$ of degree $1$, with the constraints $v^2=\te$ and $\al_1 v \al_2 = 0$ for all $\al_1,\al_2\in \Hch$. Then the algebra $\Pch[v]$ generated by $\Pch$ and all products with $v$ is endowed with a canonical differential $d$ as follows:
$$
d\al = \tilde{d}_H\al + v \al + (-1)^{|\al|}\al v\ ,\qquad dv=0\ ,
$$
where $|\al|$ is the degree of $\al\in\Pch$. One easily checks that $d^2=0$, i.e. $(\Pch[v],d)$ is a DG pro-algebra. The homomorphism $\mu_*:\Tt\Ec\to \Pch$ and the differential $d$ on $\Pch[v]$ give rise to a cocycle $\chi(\mu_*,d) \in \hom(\Omh\Tt\Ec,X(\Pch[v],d))$ defined by the same formulas as the cocycle $\chi(\rho_*,d)\in \hom(\Omh\Th\Cc,X(\Hch,d))$ of Lemma \ref{lchi}. By remark 4.4 of \cite{P10}, choose a generalized connection on $\pi:M\to B$ by fixing some horizontal distribution $H'\subset TM$. Using this generalized connection the residue morphism (\ref{resmap}) extends in an obvious fashion to a morphism
$$
\Res: \Omh\Tt\Ec\to (\Omh\Th\Cc_\top)_{[B]}
$$
On the other hand, the pullback of $H'$ on the submersion $E\times_BM\to E$ yields a compatible horizontal distribution $H''\subset T(E\times_BM)$. Choose $H=H''\cap(\id,\pi)_*^{-1}(F)$ as horizontal distribution defining the foliated connection $d_H$ (remark that the latter is generally not $\Gamma$-invariant), and thus also the DG algebra $(\Pch[v],d)$. The residue map (\ref{res}) extends to a morphism of DG-algebra $X$-complexes $\Res: X(\Pch[v],d)\to X(\Hch_\top,d)_{[E]}$ in an obvious fashion. Then a tedious computation shows that one has a commutative diagram of chain maps
$$
\xymatrix{
\Omh\Tt\Ec \ar[d]^{\Res} \ar[r]^{\chi(\mu_*,d)\quad}  & X(\Pch[v],d) \ar[d]^{\Res}  \\
(\Omh\Th\Cc_\top)_{[B]} \ar[r]^{\chi(\rho_*,d)} & X(\Hch_\top,d)_{[E]} }
$$
Replacing the differential $d$ by the superconnection $\nabla_1 = d+v$ acting by commutators on $\Pch[v]$, a JLO-type formula as (\ref{JLO1}) gives a cocycle $\chi(\mu_*,\nabla_1) \in \hom(\Omh\Tt\Ec,X(\Pch[v],d))$. By a classical transgression formula, the linear homotopy between $d$ and $\nabla_1$ shows that the cocycles $\chi(\mu_*,d)$ and $\chi(\mu_*,\nabla_1)$ are cohomologous. Now we proceed as in section \ref{sres} and enlarge the complexes $\Omh\Tt\Ec$ and $X(\Pch[v],d)$ by allowing the presence of log-polyhomogeneous pseudodifferential operators. Thus let $\Omh\Tt\Ec^1_{\log}$ and $X(\Pch[v],d)^1_{\log}$ be the complexes containing at most one power of the logarithm $\ln Q$. The cocycles $\chi(\mu_*,d)$ and $\chi(\mu_*,\nabla_1)$ extend to cohomologus cocycles in $\hom(\Omh\Tt\Ec^1_{\log}, X(\Pch[v],d)^1_{\log})$. Also the above residue morphisms extend to morphisms
\beq
\Res &:& \Omh\Tt\Ec^1_{\log} \cap \dom(\Res) \to (\Omh\Th\Cc_\top)_{[B]}\ ,\non\\ 
\Res &:& X(\Pch[v],d)^1_{\log}\cap\dom(\Res)\to X(\Hch_\top,d)_{[E]}\ ,\non
\eeq
where the domains $\dom(\Res)$ are linearly generated by differences of chains for which only the place of $\ln Q$ changes. Let $t\in [0,1]$ be a parameter, and denote by $\Om[0,1]$ the de Rham complex of differential forms over the interval, with differential $d_t$. Using the superconnection
$$
\nabla_2 = d + d_t + v + t\,\epsilon \ln Q
$$
we view the corresponding JLO cocycle $\chi(\mu_*,\nabla_2)$ in the complex $\Om[0,1]\otimes \hom(\Omh\Tt\Ec, X(\Pch[v][\epsilon],d)^1_{\log})$. Define the eta-cochain
$$
\eta(\mu_*,\nabla_2) = - \frac{\d}{\d \epsilon} \int_{t=0}^{1}\chi(\mu_*,\nabla_2)
$$
in $\hom(\Omh\Tt\Ec, X(\Pch[v],d)^1_{\log})$. The property $d_t\chi(\mu_*,\nabla_2)+ [\d,\chi(\mu_*,\nabla_2)] = 0$ implies the equality of cocycles
$$
\chi^\Res(\si_*',d_H,\ln q) = \Res\circ[\d,\eta(\mu_*,\nabla_2)]\circ \si_*
$$
in $\hom(\Omh\Th\Ac, X(\Hch_\top,d)_{[E]})$. Now if $\et=e_1\otimes\ldots\otimes e_n\in T\Ec $ is a tensor, we define its product with the left multiplier $\ln Q$ as $\ln Q\cdot\et = (\ln Q\cdot e_1)\otimes\ldots\otimes e_n$. Then define a linear map $\psi\in \hom(\Omh\Tt\Ec,\Omh\Tt\Ec^1_{\log})$ as follows:
$$
\psi(\et_0\dd \et_1 \dd \et_2\ldots \dd \et_n) = \et_0 \dd(\ln Q \cdot \et_1)\dd \et_2\ldots\dd \et_n
$$
on any $n$-form, $n\geq 1$, and $\psi(\et_0) = \ln Q\cdot \et_0$. The commutator of $\psi$ with the Hochschild boundary map $b$ on both $\Omh\Tt\Ec$ and $\Omh\Tt\Ec^1_{\log}$ is the degree -1 map
$$
[b,\psi](\et_0\dd \et_1 \dd \et_2\ldots \dd \et_n) = \et_0(\ln Q\cdot \et_1) \dd \et_2\ldots \dd \et_n - \et_0\et_1\dd (\ln Q\cdot \et_2)\ldots \dd \et_n 
$$
for $n\geq 2$, and $[b,\psi](\et_0\dd \et_1)=\et_0(\ln Q\cdot \et_1) - (\ln Q\cdot \et_0)\et_1$. Similarly $[B,\psi]$ is a simple algebraic expression involving differences of pairs of terms where only $\ln Q$ moves. Hence the range of $[\d,\psi]$, with $\d = b+B$ the total boundary of the cyclic bicomplex, is actually contained in the domain of the residue morphism. Therefore $[\d,\psi]\in \hom(\Omh\Tt\Ec,\Omh\Tt\Ec^1_{\log}\cap\dom(\Res))$ is a (non-trivial) cocycle. The diagram 
$$
\xymatrix{
\Omh\Tt\Ec \ar[d]_{\psi} \ar[dr]^{\eta(\mu_*,\nabla_2)} &  \\
\Omh\Tt\Ec^1_{\log} \ar[r]_{\chi(\mu_*,\nabla_1)\quad} & X(\Pch[v],d)^1_{\log} }
$$
is not commutative, but the only difference between $\eta(\mu_*,\nabla_2)$ and $\chi(\mu_*,\nabla_1)\circ\psi$ is that $\ln Q$ does not appear at the same places. Hence the range of the difference $\eta(\mu_*,\nabla_2) - \chi(\mu_*,\nabla_1)\circ\psi$ is contained in the domain $X(\Pch[v],d)^1_{\log}\cap\dom(\Res)$ of the residue morphism, and the cocycle $[\d,\eta(\mu_*,\nabla_2)]$ is cohomologous to the cocycle $\chi(\mu_*,\nabla_1)\circ [\d,\psi]$ in $\hom(\Omh\Tt\Ec,X(\Pch[v],d)^1_{\log}\cap\dom(\Res))$. Finally, the equality $\Res\circ\chi(\mu_*,d) = \chi(\rho_*,d)\circ\Res$ extends to an equality in $\hom(\Omh\Tt\Ec^1_{\log}\cap \dom(\Res) , X(\Hch_\top,d)_{[E]})$. Collecting everything, we get a diagram of chain maps
$$
\xymatrix{
X(\Th\Ac) \ar[r]^{\si_*\circ\gamma} & \Omh\Tt\Ec \ar[d]_{[\d,\psi]} \ar[dr]^{\quad [\d,\eta(\mu_*,\nabla_2)]} & & \\
 & \Omh\Tt\Ec^1_{\log}\cap \dom(\Res) \ar[d]^{\Res} \ar[r]_{\chi(\mu_*,d)\quad} \ar[dl] & X(\Pch[v],d)^1_{\log}\cap\dom(\Res) \ar[d]^{\Res} & \\
 X(\Th\Cc_\top)_{[B]} & (\Omh\Th\Cc_\top)_{[B]} \ar[l]_{p_X} \ar[r]^{\chi(\rho_*,d)} & X(\Hch_\top,d)_{[E]} \ar[r]^{\la'(c)} & \cc}
$$
which is commutative up to homotopy. The bottom left arrow $p_X$ is a homotopy equivalence, with inverse given by the map $\gamma: X(\Th\Cc_\top)_{[B]} \to (\Omh\Th\Cc_\top)_{[B]}$. By construction the cocycle $\varphi\in \hom(X(\Th\Cc_\top),\cc)$ representing $[N,E,\Phi,c]$ is the composition of the bottom arrows. By Theorem  \ref{tres}, the cocycle $\pi^!_G(\varphi) \in \hom(X(\Th\Ac),\cc)$ is the path surrounding the diagram via the bottom left corner, while (\ref{path}) is the upper path. \cqfd\\

\section{Dirac superconnections}\label{sdir}

Let $\La^\bullet T^*_\pi M \otimes\cc$ be the (complexified) exterior algebra over the vertical cotangent bundle associated to the submersion $\pi$, and denote by $V$ its pullback under the projection $E\times_B M\to E$. Then $V$ is a $\zz_2$-graded complex vector bundle over $E\times_B M$, and its algebra of smooth sections, which is a quotient of the algebra of all differential forms over $E\times_BM$, may be called ``vertical" differential forms. The algebra bundle of (non-compactly supported) vertical scalar symbols $\CS(E\times_B M)$ can be enlarged to a $\zz_2$-graded algebra bundle $\CS(E\times_B M,V)$ over $E$, whose fiber is the algebra of vertical pseudodifferential symbols acting on the smooth sections of $V$. Let $\PS(E\times_B M,V)$ be the subbundle of \emph{polynomial} symbols, i.e. the symbols of vertical \emph{differential} operators acting on the smooth sections of $V$. Since $\CS(E\times_B M,V)$ is an algebra bundle, there is a left representation $L:\CS(E\times_B M,V) \to \End(\CS(E\times_B M,V))$ and a right representation $R:\CS(E\times_B M,V)^{\mathrm{op}}\to \End(\CS(E\times_B M,V))$ as endomorphisms, and the two actions commute in the graded sense. In particular the graded tensor product $\CS(E\times_B M,V)\otimes \PS(E\times_B M,V)^{\mathrm{op}}$ is naturally represented in the endomorphism bundle. We  let
$$
\Lc(E\times_B M) = \im\big( \CS(E\times_B M,V)\otimes \PS(E\times_B M,V)^{\mathrm{op}}\to \End(\CS(E\times_B M,V))\big)
$$
be the range of this representation. Hence $\Lc(E\times_B M)$ is a $\zz_2$-graded algebra bundle over $E$, whose fiber is a certain algebra of linear operators acting on vertical symbols. To become familiar with these objects we introduce a local foliated coordinate system $(z^a,y^\mu)$ on an open subset $U\subset E$, such that $(y^\mu)_{\mu=1,2,\ldots}$ are the coordinates along the leaves of the foliation $F$, and $(z^a)_{a=1,2,\ldots}$ are transverse coordinates. We complete this system with vertical coordinates $(x^i)_{i=1,2,\ldots}$ on the fibers of $M$, so that $(z^a,y^\mu,x^i)$ is a local foliated coordinate system on $E\times_B M$. The vertical vectors are generated by the partial derivatives $\d/\d x^i$, and the smooth sections of the vector bundle $V$ (the vertical differential forms) are generated by products of one-forms $dx^i$ modulo horizontal forms. Denote by $\i p_i$, with $\i=\sqrt{-1}$, the Lie derivative of vertical differential forms with respect to the vector field $\d/\d x^i$. Viewing the coordinate $x^i$ as multiplication operator by the function $x^i$, these operators of even dgree fulfill the canonical commutation relations
$$
[x^i,x^j] = 0\ ,\qquad [x^i,p_j] = \i\delta^i_j\ ,\qquad [p_i,p_j] = 0\ ,
$$
and of course all commutators with the basic coordinates $z^a,y^\mu$ vanish. Let $\psi^i$ denote the operator of exterior product from the left by $dx^i$ (modulo horizontal forms) on vertical differential forms, and $\psib_i$ the operator of interior product by the vector field $\d/\d x^i$. Then $\psi,\psib$ are operators of odd degree and fulfill the canonical anticommutation relations
$$
[\psi^i,\psi^j] = 0\ ,\qquad [\psi^i,\psib_j] = \delta^i_j\ ,\qquad [\psib_i,\psib_j] = 0
$$
where the commutators are graded. Thus in local coordinates a polynomial section $a\in \cinf(E,\PS(E\times_B M,V))$ of the bundle of symbols is a smooth function of $(z,y,x,p,\psi,\psib)$ depending polynomially on $p$. Since $\psi$ and $\psib$ are odd coordinates any smooth function of them is also automatically polynomial. In the same way, a section $a\in\cinf(E,\CS(E\times_BM,V))$ of order $m$ is locally an asymptotic expansion $a\sim\sum_{j\geq 0} a_{m-j}$ where each $a_{m-j}$ is a smooth function of $(z,y,x,p,\psi,\psib)$ homogeneous of degree $m-j$ in $p$. One then shows (\cite{P9}) that the sections $s\in \cinf(E,\Lc(E\times_BM))$, which act by linear operators on $\cinf(E,\CS(E\times_BM,V))$, are asymptotic expansions of the partial derivatives $\d/\d x$ and $\d/\d p$ of the form
\be 
s = \sum_{|\al|=0}^k \sum_{|\beta|=0}^{\infty} \sum_{|\gamma|\geq 0} \sum_{|\delta|\geq 0} (s_{\al,\beta,\gamma,\delta})_L (\psi)_R^\gamma(\psib)^\delta_R \Big(\frac{\d}{\d x}\Big)^\al \Big(\frac{\d}{\d p}\Big)^\beta  \label{asym}
\ee
where $\al,\beta,\gamma,\delta$ are multi-indices, and $s_{\al,\beta,\gamma,\delta}$ is a local section of $\CS(E\times_BM,V)$. Let $\eps$ be an indeterminate (of even parity, not to be confused with the previous odd parameter $\epsilon$) and consider the $\zz_2$-graded algebra bundle of formal power series $\Sc(E\times_BM) = \Lc(E\times_BM)[[\eps]]$ in $\eps$. One defines a filtered, $\zz_2$-graded algebra sub-bundle
\be 
\Dc(E\times_BM)=\bigcup_{m\in\rr} \Dc^m(E\times_BM) \subset \Sc(E\times_BM)
\ee
as follows: a formal series $s=\sum_{k\geq 0} s_k\eps^k$ is a section of $\Dc^m(E\times_BM)$ if each coefficient $s_k$ is locally an asymptotic expansion (\ref{asym}), where the symbol $s_{\al,\beta,\gamma,\delta}$ has order $\leq m +(k+|\beta|-3|\al|)/2$. According to this filtration, the partial derivatives $\d/\d x$ have degree $m=3/2$, the partial derivatives $\d/\d p$ and $\eps$ both have degree $m=-1/2$, and a symbol $a\in\cinf(E,\CS^m(E\times_BM,V))$ in the \emph{left} representation $a_L$ has degree $m$. Following \cite{P9} Definition 5.1, a \emph{generalized Dirac operator} as an odd section $D\in \cinf(E,\Dc(E\times_B M))$ which in local coordinates reads
$$
D = \i\eps(\psi^i)_R\Big(\frac{\d}{\d x^i} + \ldots\Big) + (\psib_i)_R\Big(\frac{\d}{\d p_i}+\ldots\Big)
$$
where the dots are operators of lower degree (according to the filtration of $\Dc(E\times_BM)$) given by expansions in powers of the partial derivatives $\d/\d p$. Summation over the repeated indices $i$ is understood. Using a partition of unity one shows that such operators always exist globally on $E\times_BM$. More importantly, the general form of a Dirac operator above \emph{is preserved under any change of local coordinates} compatible with the submersion. This makes the present formalism particularly well-adapted to groupoid actions. The square of $D$ is a generalized Laplacian taking locally the form
$$
\Delta = -D^2 = \i\eps \frac{\d}{\d x^i}\frac{\d}{\d p_i} + \ldots
$$
Due to the presence of an overall factor $\eps$ in the Laplacian, the heat operator $\exp(- D^2) \in \cinf(E,\Sc(E\times_BM))$ is a well-defined formal power series. Moreover, a Duhamel-like expansion holds for the heat operator of perturbed Laplacians \cite{P9}. Then let $\Tc(E\times_BM)$ be the vector subbundle of $\Sc(E\times_BM)$ whose sections are of the form $s\exp(-D^2)$, for all sections $s\in\cinf(E,\Dc_c(E\times_BM))$ with \emph{compact} vertical support. One shows as in \cite{P9} that $\Tc(E\times_BM)$ is a $\zz_2$-graded $\Dc(E\times_BM)$-bimodule, and that there exists a \emph{canonical} graded trace
\be 
\Tr_s: \cinf(E,\Tc(E\times_BM)) \to \cinf(E) \label{trs}
\ee
coming from the fiberwise Wodzicki residue. $\Tc(E\times_BM)$ is called the bimodule of trace-class operators. \\
In local coordinates, the horizontal distribution $H$ associated to a choice of decomposition $(\id,\pi)_*^{-1}(F) = H \oplus \ker(\id,\pi)_*$ is the intersection of the kernels of the collection of $1$-forms $dx^i - \om^i_\mu dy^\mu$ (summation over repeated indices), where $\om^i_\mu$ are scalar functions over $E\times_BM$. The associated connection $d_H$ on the bundle of vertical scalar functions is locally expressed by
$$
d_H =  dy^\mu \Big(\frac{\d}{\d y^\mu} +  \i\om^i_\mu p_i\Big) \quad \mbox{on}\quad \cinf(E,\La^\bullet F^*\otimes\cinfc(E\times_BM))\ .
$$
The curvature of $d_H$ is the horizontal 2-form $\te = \frac{1}{2} dy^\mu \wedge dy^\nu \te_{\mu\nu}^i \frac{\d}{\d x^i}$ with values in vertical vector fields, whose components read 
$$
\te^i_{\mu\nu} = \frac{\d\om^i_\nu}{\d y^\mu} + \om^j_\mu\frac{\d\om^i_\nu}{\d x^j}\ .
$$
Now we promote $d_H$ to a derivation on \emph{all} vertical differential forms. The new local expression has to be modified as follows: 
$$
d_H = dy^\mu \frac{\d}{\d y^\mu} + dy^\mu \Big(\i\om^i_\mu p_i+\frac{\d \om^i_\mu}{\d x^j}\psi^j\psib_i\Big) - \frac{1}{2} dy^\mu \wedge dy^\nu (\te_{\mu\nu}^i\psib_i)\ .
$$
Note that $\i\om^i_\mu p_i+\frac{\d \om^i_\mu}{\d x^j}\psi^j\psib_i$ is the Lie derivative of vertical differential forms with respect to the vector field $\om^i_\mu\frac{\d}{\d x^i}$, while $\te_{\mu\nu}^i\psib_i$ is the interior product by the vector field $\te^i_{\mu\nu}\frac{\d}{\d x^i}$. Both are smooth sections of the bundle $\PS^1(E\times_BM,V)$ of vertical polynomial symbols (i.e. differential operators) locally defined over $U$. The action by commutator $\tilde{d}_H=[d_H,\ ]$ is an odd derivation on the space of sections $\cinf(E,\La^\bullet F\otimes\CS(E\times_BM,V))$, explicitly 
\beq
\tilde{d}_H &=& dy^\mu \frac{\d}{\d y^\mu} + dy^\mu \Big(\i\om^i_\mu p_i+\frac{\d \om^i_\mu}{\d x^j}\psi^j\psib_i\Big)_L - dy^\mu \Big(\i\om^i_\mu p_i+\frac{\d \om^i_\mu}{\d x^j}\psi^j\psib_i\Big)_R \non\\
&& \qquad\qquad - \frac{1}{2}dy^\mu \wedge dy^\nu (\te_{\mu\nu}^i\psib_i)_L + \frac{1}{2}dy^\mu \wedge dy^\nu (\te_{\mu\nu}^i\psib_i)_R\ . \non
\eeq
Except for the derivative term $dy^\mu \frac{\d}{\d y^\mu}$, all other terms are local sections of $\La^\bullet F\otimes\Dc(E\times_BM)$ over $U$. Hence using the formalism of superconnections we can add genuine global sections of this algebra bundle to $\tilde{d}_H$. We consider the family of \emph{Dirac superconnections}
\be 
\DD = \i\eps (\tilde{d}_H + A) + D\ ,
\ee
where $A\in \cinf(E,\La^1 F\otimes\Dc(E\times_BM))$, and $D\in  \cinf(E,\Dc(E\times_BM))$ is a generalized Dirac operator. Since the complementary part of $D$ has an horizontal form degree $\geq 1$, a Duhamel expansion shows that the heat operator $\exp(-\DD^2)\in \cinf(E,\La^\bullet F\otimes\Tc(E\times_BM))$ is a trace-class section. Now observe that the bundles $\Sc(E\times_BM)$, $\Dc(E\times_BM)$, $\Tc(E\times_BM)$, etc... over $E$ are all $\Gamma$-bundles. In particular we can form the convolution algebra 
\be
\Uc = \cinf_p(E, r^*(\La^\bullet F^*\otimes\Dc_c(E\times_B M))) \cp\Gamma \ .
\ee
It naturally inherits $\zz_2$-graduation from $\La^\bullet F$ and $\Dc_c(E\times_BM)$.  $\Uc$ is a bimodule over the algebra of sections $\cinf(E,\La^\bullet F^*\otimes\Dc_c(E\times_B M))$, and if the horizontal distribution $H$ is $\Gamma$-invariant, the commutator $[\tilde{d}_H,\ ]$ defines an odd derivation on $\Uc$. Following the general recipe we construct an extension of the convolution algebra. Define the vector space $\Vc = \bigoplus_{n\geq 1} \Vc_n$, where
$$
\Vc_n = \cinf_p(\Gamma^{(n)}, r^*_1(\La^\bullet F^*\otimes\Dc_c(E\times_B M)) \otimes r^*_1|\La^{\max} A^*\Gamma|\otimes \ldots \otimes r^*_n|\La^{\max} A^*\Gamma|)
$$
As a vector space $\Vc_1$ is isomorphic to $\Uc$. A product $\Vc_{n_1}\times\Vc_{n_2} \to \Vc_{n_1+n_2}$ is defined as usual, as well as the multiplication map $\Vc\to\Uc$. Hence $\Vc$ is a $\zz_2$-graded extension of $\Uc$. The commutator $[\tilde{d}_H,\ ]$ lifts uniquely to an odd derivation on $\Vc$. More generally the commutator with any Dirac superconnection yields an odd derivation, where $\cinf(E,\La^\bullet F^*\otimes\Dc_c(E\times_B M))$ acts by multipliers on $\Vc$ in the obvious way. Finally we denote by $\Wc = \ker(\Vc\to\Uc)$ the kernel of the multiplication map, which is stable by $[\tilde{d}_H,\ ]$, and by $\Vch$ the $\Wc$-adic completion of $\Vc$. \\
Since the space of vertical 0-forms on $E\times_BM$ is a direct summand in the space of all vertical forms, the algebra bundle of scalar symbols $\CS(E\times_BM)$ sits naturally as an algebra subbundle of $\CS(E\times_BM,V)$. The latter may further be identified with a subbundle of $\Dc(E\times_BM)$ through the left representation $L$. Hence one gets a canonical inclucion $\CS_c(E\times_BM) \hookrightarrow \Dc_c(E\times_BM)$, which in turn induces an homomorphism of pro-algebras $\Pch\hookrightarrow \Vch$. Composing with the homomorphism $\si_*':\Th\Ac\to\Pch$ constructed above one gets a representation
\be 
\si_*''=L\circ\si_*'\ :\ \Th\Ac \to \Pch \hookrightarrow \Vch\ .
\ee
Choose as above an elliptic section $Q\in \cinf(E,\CL^1(E\times_BM))$, with symbol $q\in \cinf(E,\CS^1(E\times_BM))$. Extend $q$ to an elliptic symbol acting on vertical differential forms $\tilde{q} \in \cinf(E,\CS^1(E\times_BM,V))$, requiring that the leading symbol of $\tilde{q}$ remains of scalar type. Let $\epsilon$, $\epsilon^2=0$ be the odd parameter introduced above. For any Dirac superconnection $\DD$, the new superconnection
$$
\nabla = \DD + \epsilon\ln\tilde{q}_L
$$
acting on the algebra $\Vch[\epsilon]$ by graded commutators is a graded derivation. We use the homomorphism $\si_*''$ and the superconnection $\nabla$ to construct a cocycle $\chi^{\Tr_s}(\si_*'',\nabla) \in \hom(\Omh\Th\Ac, X(\Hc_\top[\epsilon],d)_{[E]})$ by a JLO-type formula. We first extend the graded trace (\ref{trs}) to an $X$-complex map
$$
\Tr_s\ :\  X(\Vch[\epsilon]) \cap \dom (\Tr_s) \to X(\Hch_\top[\epsilon],d)_{[E]}\ ,
$$
taking into account a rescaling factor of $(\i\eps)^{-1}$ for each leafwise form degree in $\Hc$. Then on any $n$-form $\ah_0\dd\ah_1\ldots\dd\ah_n$ set
\beq
\lefteqn{\chi^{\Tr_s}(\si_*'',\nabla)(\ah_0\dd\ah_1\ldots\dd\ah_n)=}\\
&&\hspace{-1cm}\sum_{i=0}^n (-)^{i(n-i)} \int_{\Delta_{n+1}} \Tr_s\big(e^{-t_{i+1}\nabla^2} [\nabla,\si_{i+1}''] \ldots  e^{-t_{n+1}\nabla^2} \si_0'' \,e^{-t_0\nabla^2} [\nabla,\si_1''] \ldots e^{-t_i\nabla^2}\big) dt \non\\
&+&  \sum_{i=1}^n \int_{\Delta_n} \Tr_s\big(\nat \si_0''\, e^{-t_0\nabla^2} [\nabla,\si_1''] \ldots  e^{-t_{i-1}\nabla^2} \dd \si_i'' \,e^{-t_i\nabla^2} \ldots [\nabla,\si_n''] e^{-t_n\nabla^2}\big)dt \non
\eeq
where $\si_i'' = \si_*''(\ah_i)\in \Vch$ for all $i$. As above the $\epsilon$-component of $\chi^{\Tr_s}(\si_*'',\nabla)$ 
\be 
\chi^{\Tr_s}(\si_*'',\DD,\ln\tilde{q}_L) = \frac{\d}{\d \epsilon} \chi^{\Tr_s}(\si_*'',\DD+\epsilon\ln\tilde{q}_L)
\ee
is a cocycle in $\hom(\Omh\Th\Ac,X(\Hc_\top,d)_{[E]})$, whose cohomology class does not depend on the choice of Dirac superconnection $\DD$ and elliptic symbol $\tilde{q}$. Choosing genuine Dirac superconnections thus allows to build cohomologous cocycles. We now exploit this fact. Let $d_V$ be the vertical part of the de Rham operator on $E\times_BM$, acting on the vertical differential forms. One thus has $d_V\in \cinf(E,\PS(E\times_MB,V))$ and in local coordinates
$$
d_V = \i p_i\psi^i\ .
$$
Of course $d_V$ is completely canonical and $d_V^2=0$. The choice of horizontal distribution $H$ allows to identify $V$ with a subbundle of $\La^\bullet T^*(E\times_BM)$, and according to this identification the total de Rham differential on $E\times_BM$ is exactly the sum $d_H+d_V$. Thus in particular $(d_H+d_V)^2 = d_H^2+[d_H,d_V]=0$. This can be explicitly checked in local coordinates using the formulas above. Next, taking the image of $d_V$ under the right representation $R$ yields a global section of $\Dc(E\times_BM)$. In local coordinates one has
$$
(\i p_i\psi^i)_R = -(\psi^i)_R\Big(\frac{\d}{\d x^i} - (\i p_i)_L\Big)
$$
using the commutation relations. A Dirac operator $D$ is called of \emph{de Rham-Dirac} type \cite{P9} if it is exactly given by
$$
D = \i\eps(\psi^i)_R\Big(\frac{\d}{\d x^i} - (\i p_i)_L\Big) + (\psib_i)_R\Big(\frac{\d}{\d p_i}+\ldots\Big)
$$
where the dots represent an expansion in higher powers of $\frac{\d}{\d p}$ ensuring that $D$ is globally defined. Since $d_V$ is completely canonical, the term proportional to $\psi_R$ is always $\Gamma$-invariant. Only the term proportional to $\psib_R$ may not be $\Gamma$-invariant, however using the properness of the action we can always find a $\Gamma$-invariant de Rham-Dirac operator.

\begin{proposition}\label{pderham}
Let $\DD = \i\eps \tilde{d}_H + D$ be a Dirac superconnection constructed from a $\Gamma$-invariant horizontal distribution $H$ and a $\Gamma$-invariant de Rham-Dirac operator $D$. Let $q$, $\tilde{q}$ be $\Gamma$-invariant elliptic symbols of scalar type. Then one has the equality of cocycles in $\hom(\Omh\Th\Ac,X(\Hc_\top,d)_{[E]})$
\be 
\chi^{\Tr_s}(\si_*'',\DD,\ln \tilde{q}_L) = \chi^{\Res}(\si_*',d_H,\ln q)\ .
\ee
\end{proposition}
{\it Proof:} Let $\Pi \in \cinf(E,\PS(E\times_B,V))$ be the natural projection operator of vertical differential forms onto their 0-degree component (scalar functions). Write $D=-\i\eps (d_V)_R + \nablab$ where $\nablab$ contains all the terms proportional to $\psib_R$. Since left and right representations commute one has $[\tilde{d}_H,\si_*''(\ah)] = ([d_H,\si_*'(\ah)\Pi])_L$ and $[D,\si_*''(\ah)] = [\nablab,\si_*''(\ah)]$ for any $\ah\in\Th\Ac$. Thus
$$
[\DD,\si_*''(\ah) ] = ([d_H,\si_*'(\ah)\Pi])_L + [\nablab,\si_*''(\ah)]\ .
$$
On the other hand $\DD^2 = -\eps^2(\tilde{d}_H)^2 + \i\eps [\tilde{d}_H,D] + D^2$, with $(\tilde{d}_H)^2 = \te_L - \te_R$ and $\i\eps [\tilde{d}_H,D] = \eps^2([d_H,d_V])_R + \i\eps [\tilde{d}_H,\nablab] =-\eps^2 \te_R + \i\eps [\tilde{d}_H,\nablab]$. Hence the Laplacian reads
$$
-\DD^2 = \eps^2 \te_L - D^2 - \i\eps [\tilde{d}_H,\nablab]
$$
At this point we can apply verbatim the proof of \cite{P9} Proposition 6.4 (and 6.7), showing that the terms $[\nablab,\si_*''(\ah)]$ and $[\tilde{d}_H,\nablab]$ do not contribute, and that the effect of the graded trace $\Tr_s$ is to reduce the JLO cocycle $\chi(\si_*'',\DD,\ln \tilde{q}_L)$ to the JLO cocycle $\chi(\si_*',d_H,\ln q)$. \cqfd\\

Choose a torsion-free affine connection on the manifold $E\times_BM$ and take its restriction $\nabla^T$ to the vertical tangent bundle $E\times_BT_\pi M$. In local coordinates we can write $\nabla^T = dz^a\nabla_a + dy^\mu\nabla_\mu + dx^i\nabla_i$, and because we are interested to differentiate in the directions of the leaves ($z$ fixed) we only retain the covariant derivatives $\nabla^T_\mu$ and $\nabla^T_i$ in the horizontal and vertial directions respectively. Their effect on vertial vector fields is expressed in terms of the Christoffel symbols $\Gamma^k_{i\mu}=\Gamma^k_{\mu i}$ and $\Gamma^k_{ij}=\Gamma^k_{ji}$ of the affine connection:
$$
\nabla^T_\mu\Big(\frac{\d}{\d x^j}\Big) = \Gamma^k_{\mu j}\frac{\d}{\d x^k}\ , \qquad \nabla^T_i\Big(\frac{\d}{\d x^j}\Big) = \Gamma^k_{ij}\frac{\d}{\d x^k}\ .
$$
In particular we denote $\nabla^T_i\om^k = \frac{\d \om^k}{\d x^i} + \Gamma^k_{ij}\om^j$ the covariant derivative of the 1-form $\om^k = \om^k_\mu dy^\mu$ defining the horizontal distribution $H$. In local coordinates the curvature of the affine connection acts on vertical vector fields by
$$
\big(\nabla^T\big)^2\Big(\frac{\d}{\d x^l}\Big) = R^k_l\frac{\d}{\d x^k}\ ,
$$
where the coefficients $R^k_l$ are leafwise 2-forms over $E\times_BM$. We decompose them into purely horizontal, mixed and purely vertical components  
$$
R^k_l = \frac{1}{2} R^k_{l\mu\nu}dy^\mu\wedge dy^\nu + R^k_{l\mu j}dy^\mu\wedge dx^j + \frac{1}{2} R^k_{lij}dx^i\wedge dx^j\ ,
$$
where the components of the curvature tensor are expressed as usual via the Christoffel symbols
\beq
R^k_{l\mu\nu} &=& [\nabla^T_\mu,\nabla^T_\nu]^k_l \ =\ \frac{\d \Gamma^k_{\nu l}}{\d y^\mu} - \frac{\d \Gamma^k_{\nu l}}{\d y^\mu} + \Gamma^k_{\mu m}\Gamma^m_{\nu l} - \Gamma^k_{\nu m}\Gamma^m_{\mu l} \non\\
R^k_{l\mu j} &=& [\nabla^T_\mu,\nabla^T_j]^k_l \ =\ \frac{\d \Gamma^k_{jl}}{\d y^\mu} - \frac{\d \Gamma^k_{\mu l}}{\d x^j} + \Gamma^k_{\mu m}\Gamma^m_{jl} - \Gamma^k_{jm}\Gamma^m_{\mu l} \non\\
R^k_{lij} &=& [\nabla^T_i,\nabla^T_j]^k_l \ =\ \frac{\d \Gamma^k_{jl}}{\d x^i} - \frac{\d \Gamma^k_{il}}{\d x^j} + \Gamma^k_{im}\Gamma^m_{jl} - \Gamma^k_{jm}\Gamma^m_{il}\ . \non
\eeq
\begin{definition}
A Dirac superconnection $\DD = \i\eps (\tilde{d}_H + A) + D$ is called \emph{affiliated to the affine connection $\nabla^T$ and the horizontal distribution $H$} if in local coordinates
\beq
A &=& (\nabla^T_j\om^k + \Gamma^k_{\mu j}dy^\mu )_L \Big((p_k)_L\frac{\d}{\d p_j} + (\psib_k\psi^j)_L+\ldots \Big) \non\\
&& + (\psi^i)_R \big(\psib_k (\nabla^T_i\om^k + \Gamma^k_{i\mu}dy^\mu ) +\ldots \big)_L \non\\
D &=& \i\eps (\psi^i)_R \Big(\frac{\d}{\d x^i} +(\Gamma^k_{ij}p_k)_L\frac{\d}{\d p_j} + (\Gamma^k_{ij}\psib_k\psi^j)_L+\ldots \Big) + (\psib_i)_R\Big(\frac{\d}{\d p_i}+\ldots \Big)\non
\eeq
where the dots denote an expansion in higher powers of $\frac{\d}{\d p}$. 
\end{definition}
The global existence of such operators $A$ and $D$ on $E\times_BM$ is proved as usual by gluing together local operators by means of a partition of unity. From now on let us abusively denote by $\pi$ the submersion $E\times_B S^*_\pi M\to E$, and by $\pi_*$ its tangent map. The preimage $\pi_*^{-1}F$ of the subbundle $F\subset TE$ under the tangent map defines a foliation on $E\times_B S^*_\pi M$. Pulling back the vertical tangent bundle and its affine connection from $E\times_BM$ to $E\times_B S^*_\pi M$, we can view the above curvature $R$ as a leafwise 2-form on $E\times_B S^*_\pi M$ with values in the endomorphisms of the vertical tangent bundle. By Chern-Weil theory, the Todd class of the complexified vertical tangent bundle $\Td((E\times_BT_\pi M)\otimes\cc)$ is represented by the closed leafwise differential form
\be 
\Td(\i R/2\pi) = \det\left(\frac{\i R/2\pi}{e^{\i R/2\pi} -1}\right)\ \in \cinf(E\times_B S^*_\pi M,\La^{4\bullet}\pi_*^{-1}F)\ ,\label{todd}
\ee
which is a polynomial in the Pontryagin classes. Next, the submersion $\pi:E\times_B S^*_\pi M\to E$ being $\Gamma$-equivariant, it extends to a morphism of Lie groupoids $\pi: (E\times_B S^*_\pi M)\cp \Gamma\to \Gamma$. Its tangent map $\pi_*$ is a morphism of the corresponding tangent groupoids. Hence the preimage of the flat connection $\Phi\subset \Gamma$,
\be 
\pi_*^{-1}(\Phi) \subset T((E\times_B S^*_\pi M)\cp\Gamma)\ ,
\ee
is a flat connection on the groupoid $(E\times_B S^*_\pi M)\cp\Gamma$, and the latter acts on leafwise differential forms over $E\times_BS^*_\pi M$. If we start from a $\Gamma$-equivariant affine connection $\nabla^T$, which is always possible by the properness of the action, then the Todd form $\Td(\i R/2\pi)\in \cinf(E\times_B S^*_\pi M,\La^{4\bullet}\pi_*^{-1}F)$ is invariant, and a classical homotopy argument shows that its class in the cohomology of invariant leafwise differential forms is independent of the choice of connection. Also note that one can further assume $\nabla^T$ to be a $\Gamma$-equivariant \emph{metric} connection, in the sense that it preserves an invariant scalar product on the vertical tangent bundle. Now define the convolution algebra
\be
\Xc = \cinf_p(E\times_B S^*_\pi M ,\La^\bullet\pi_*^{-1}F) \cp \Gamma \ .
\ee
It comes equipped with the grading of leafwise differential forms, and the leafwise de Rham operator $d$ which squares to zero. Hence $(\Xc,d)$ is a DG algebra. As usual we build an extension $\Yc = \bigoplus_{n\geq 1} \Yc_n$ with 
$$
\Yc_n = \cinf_p(((E\times_B S^*_\pi M)\cp\Gamma)^{(n)}, r^*_1\La^\bullet \pi_*^{-1}F^* \otimes r^*_1|\La^{\max} A^*\Gamma|\otimes \ldots \otimes r^*_n|\La^{\max} A^*\Gamma|)
$$
where for notational simplicity we identify the fiber of the algebroid $A\Gamma$ with that of $A((E\times_B S^*_\pi M)\cp\Gamma)$. Then $\Yc_1=\Xc$ as a vector space, there is a product $\Yc_{n_1}\times\Yc_{n_2}\to\Yc_{n_1+n_2}$ and a multiplication map $\Yc\to\Xc$. Moreover the differential $d$ on $\Xc$ extends uniquely to a differential on $\Yc$. We let $\Zc=\ker( \Yc\to\Xc)$ be the kernel of the multiplication map which is a DG ideal in $\Yc$, and denote by $\Ych$ the $\Zc$-adic completion of $\Yc$. The DG pro-algebra $(\Ych,d)$ is the analogue, for the groupoid $(E\times_B S^*_\pi M)\cp\Gamma$ with connection $\pi_*^{-1}(\Phi)$, of the DG pro-algebra $(\Hch,d)$ introduced in section \ref{sgeo} for the groupoid $\Gamma$ with connection $\Phi$. In order to apply Lemma \ref{lchi} in this context, we need to build an homomorphism from $\Th\Ac$ to $\Ych$. To this end, observe that the cut-off function previously chosen on the submersion $\nu:N\to B$ yields, by pullback to the cosphere bundle, a cut-off function on the submersion $N\times_BS^*_\pi M \to S^*_\pi M$, whence an homomorphism of algebras
$$
\Ac = \cinfc(S^*_\pi M)\cp G \to \cinfc(N\times_B S^*_\pi M)\cp \nu^*G
$$ 
realizing the Morita equivalence between the corresponding groupoids. Then as before the $\Gamma$-equivariant map $\eta:E\to N$ induces a pullback homomorphim $\cinfc(N\times_BS^*_\pi M)\cp \nu^*G \to \cinf_p(E\times_BS^*_\pi M)\cp \Gamma\subset \Xc$. The resulting homomorphism $\rho:\Ac \to \Xc$ extends as usual to an homomorphism of pro-algebras
\be 
\rho_* : \Th\Ac \to \Ych \label{homo}
\ee
Remark that $\rho:\Ac\to\Xc$ could as well be obtained by composition of the linear map $\si:\Ac\to\Ec$ with the homomorphism $\mu:\Ec\to \Oc^0$ (where the range only contains formal symbols of order $\leq 0$), followed by the projection homomorphism $\Oc^0\to \Xc$ onto leading symbols (i.e. scalar functions over the cosphere bundle). Hence (\ref{homo}) is the composition of the homomorphism $\si_*':\Th\Ac \to \Pch^0$ with the projection $\Pch^0\to \Ych$ onto leading symbols. Lemma \ref{lchi} yields a chain map $\chi(\rho_*,d) \in \hom(\Omh\Th\Ac,X(\Ych,d))$, defined on any $n$-form $\ah_0\dd\ah_1\ldots\dd\ah_n\in\Om^n\Th\Ac$ by 
\beq
\lefteqn{\chi(\rho_*,d)(\hat{a}_0 \dd\hat{a}_1\ldots \dd\hat{a}_n)=} \non\\
&& \frac{1}{(n+1)!} \sum_{i=0}^n (-1)^{i(n-i)} d\rho_*(\hat{a}_{i+1}) \ldots d\rho_*(\hat{a}_n) \, \rho_*(\hat{a}_0) \, d\rho_*(\hat{a}_1) \ldots d\rho_*(\hat{a}_i) \non \\
&& + \frac{1}{n!} \sum_{i=1}^n \nat\big(\rho_*(\hat{a}_0) \, d\rho_*(\hat{a}_1) \ldots \dd \rho_*(\hat{a}_i) \ldots d\rho_*(\hat{a}_n) \big)\ . \non
\eeq
Also note that the wedge product by a closed, invariant, leafwise differential form on $E\times_BS^*_\pi M$ defines in an obvious way an endomorphism of the localized complex $X(\Ych_\top,d)_{[E\times_BS^*_\pi M]}$. Moreover, integration of differential forms along the fibers of the submersion $E\times_BS^*_\pi M\to E$ gives rise to a chain map of total complexes
$$
\int_{S^*_\pi M} \ :\ X(\Ych_\top,d)_{[E\times_BS^*_\pi M]} \to X(\Hch_\top,d)_{[E]}\ .
$$

\begin{proposition}\label{paffil}
Let $\DD = \i\eps (\tilde{d}_H + A) + D$ be a Dirac superconnection affiliated to a $\Gamma$-equivariant metric connection $\nabla^T$ on the vertical tangent bundle and a $\Gamma$-invariant horizontal distribution $H$. Let $\tilde{q}$ be a $\Gamma$-invariant elliptic symbol of scalar type. Then one has the equality of cocycles in $\hom(\Omh\Th\Ac,X(\Hc_\top,d)_{[E]})$
\be 
\chi^{\Tr_s}(\si_*'',\DD,\ln \tilde{q}_L) = \int_{S^*_\pi M}\Td(\i R/2\pi) \wedge\chi(\rho_*,d)\ ,
\ee
where $R$ is the curvature $2$-form of $\nabla^T$.
\end{proposition}
{\it Proof:} This is a direct generalization of \cite{P9} Theorem 6.5, stating that for Dirac operators affiliated to metric connections the JLO cocycle $\chi^{\Tr_s}(\si_*'',\DD,\ln \tilde{q}_L)$ only involves the \emph{leading symbols} of its arguments, so has a simple expression as an integral of ordinary differential forms over the cosphere bundle. Indeed one has
\beq
\lefteqn{[\DD,\si_*''(\ah)] = \i\eps dy^\mu \Big( \frac{\d \si_*'(\ah)}{\d y^\mu}\Pi + \Gamma^k_{\mu j}p_k\frac{\d\si_*'(\ah)}{\d p_j}\Pi - \Pi\psib_k\nabla^T_\mu \om^k\si_*'(\ah) \Big)_L } \non\\
&+& \i\eps (\psi^i+\om^i)_R \Big( \frac{\d \si_*'(\ah)}{\d x^i}\Pi + \Gamma^k_{ij} p_k\frac{\d\si_*'(\ah)}{\d p_j}\Pi - \Pi\psib_k(\nabla^T_i \om^k + \Gamma^k_{i\mu}dy^\mu)\si_*'(\ah) \Big)_L \non\\
&+& (\psib_i)_R \Big(\frac{\d \si_*'(\ah)}{\d p_i}\Pi\Big)_L + \ldots \non
\eeq
for any $\ah\in \Th\Ac$, where the dots denote expansions in higher powers of $\d/\d p$. On the other hand the Laplacian reads
\beq
\lefteqn{-\DD^2 = \i\eps \Big(\frac{\d}{\d x^i}\frac{\d}{\d p_i} + (\Gamma^k_{ij})_L ((\psi^i+\om^i)\psib_k)_R \frac{\d}{\d p_j}  + (\Gamma^k_{\mu j})_L dy^\mu(\psib_k)_R \frac{\d}{\d p_j} +\ldots \Big) } \non\\
&+& \frac{\eps^2}{2} dy^\mu \wedge dy^\nu \Big((R^k_{l\mu\nu}p_k)_L \frac{\d}{\d p_l} + (R^k_{l\mu\nu}\psib_k)_L ((\psi^l+\om^l)_L-(\psi^l+\om^l)_R) \Big) \non\\
&+& \eps^2 dy^\mu (\psi^j+\om^j)_R \Big((R^k_{l\mu j}p_k)_L \frac{\d}{\d p_l} + (R^k_{l\mu j}\psib_k)_L ((\psi^l+\om^l)_L-(\psi^l+\om^l)_R) \Big) \non\\
&+& \frac{\eps^2}{2} (\psi^i+\om^i)_R(\psi^j+\om^j)_R \Big((R^k_{lij}p_k)_L \frac{\d}{\d p_l} + (R^k_{lij}\psib_k(\psi^l+\om^l))_L \Big) + \ldots \non
\eeq
The identities $\psib_k\psi^l\Pi=\delta^k_l\Pi = \Pi\psib_k\psi^l$ and $\psib_k\Pi = 0 = \Pi\psi^l$ for all indices $k,l$, together with the fact that $\nabla^T$ is a metric connection ($R^k_{k\mu\nu}=R^k_{k\mu j} = R^k_{kij}=0$), allows to show as in \cite{P9} Theorem 6.5 that the terms involving $\psib_L$ do not contribute. The rest of the proof follows exactly the lines of \cite{P9} Theorems 6.5 and 6.8. \cqfd\\

Now remark that the bicomplex of equivariant leafwise differential forms $C^\bullet((E\times_B S^*_\pi M)\cp\Gamma,\La^\bullet\pi_*^{-1}F)$ is naturally a module (for the wedge product) over the algebra of closed invariant leafwise differential forms. In particular, if $c\in C^\bullet(\Gamma,\La^\bullet F)$ is a normalized total cocycle, we can take its pullback $\pi^*(c)$ under the projection $\pi:E\times_B S^*_\pi M\to E$, and the product
\be 
\Td(\i R/2\pi)\wedge \pi^*(c)\ \in C^\bullet((E\times_B S^*_\pi M)\cp\Gamma,\La^\bullet\pi_*^{-1}F)
\ee
is a normalized total cocycle. Clearly the latter descends to a cup-product on the corresponding cohomologies. Collecting all the preceding results we obtain

\begin{theorem}\label{tgeo}
Let $G\rightrightarrows B$ be a Lie groupoid acting on a surjective submersion $\pi:M\to B$. The excision map localized at units
$$
\pi^!_G\ :\ HP_\top^{\bullet}(\cinfc(B)\cp G)_{[B]} \to HP_\top^{\bullet+1}(\cinfc(S^*_\pi M\rtimes G))_{[S^*_\pi M]}
$$
sends the cyclic cohomology class of a \emph{proper} geometric cocycle $(N,E,\Phi,c)$ to the cyclic cohomology class
$$
\pi^!_G([N,E,\Phi,c]) = [N\times_B S^*_\pi M \, , \, E\times_B S^*_\pi M \, , \,  \pi_*^{-1}(\Phi) \, , \, \Td(T_\pi M\otimes\cc) \wedge \pi^*(c)] 
$$
where $\Td(T_\pi M\otimes\cc)$ is the Todd class of the complexified vertical tangent bundle in the invariant leafwise cohomology of $E\times_B S^*_\pi M$.
\end{theorem}
{\it Proof:} By Proposition \ref{pcompo} the class $\pi^!_G([N,E,\Phi,c])$ is represented by the cocycle $\la'(c)\circ \chi^{\Res}(\si_*',d_H,\ln q)\circ\gamma$ in $\hom(X(\Th\Ac,\cc)$. Then Proposition \ref{pderham} implies that this cocycle is cohomologous to $\la'(c)\circ\chi^{\Tr_s}(\si_*'',\DD,\ln \tilde{q}_L)\circ\gamma$ for any choice of Dirac superconnection $\DD$. By Proposition \ref{paffil} it is also cohomolgous to $\la'(\Td(T_\pi M\otimes\cc) \wedge \pi^*(c)) \circ \chi(\rho_*,d)\circ \gamma$, which is precisely the construction of the class $[N\times_B S^*_\pi M \, , \, E\times_B S^*_\pi M \, , \,  \pi_*^{-1}(\Phi) \, , \, \Td(T_\pi M\otimes\cc) \wedge \pi^*(c)]$.   \hfill\cqfd\\

Combining Theorems \ref{tres} and \ref{tgeo} yields a commutative diagram computing explicitly the excision map of the fundamental pseudodifferential extension on the range of proper geometric cocycles localized at units:
\be 
\vcenter{\xymatrix{
HP^\bullet(\cinfc(B,\CL_{c}^{-1}(M))\cp G) \ar[r]^{\quad E^*} & HP^{\bullet+1}(\cinfc(S^*_\pi M) \rtimes G)  \\
HP_\top^{\bullet}(\cinfc(B)\cp G)_{[B]} \ar[u]^{\tau_*} \ar[r]^{\pi^!_G\qquad} & HP_\top^{\bullet+1}(\cinfc(S^*_\pi M) \rtimes G)_{[S^*_\pi M]} \ar[u] }}
\ee
This together with the adjointness theorem in $K$-theory \cite{Ni,P8} gives the following
\begin{corollary}
Let $P\in M_{\infty}(\CL_{c}^0(M)\cp G)^+$ be an elliptic operator and let $(N,E,\Phi,c)$ be a proper geometric cocycle localized at units for $G$. The $K$-theoretical index $\Ind([P])\in K_0(\CL_{c}^{-1}(M)\cp G)$ evaluated on the cyclic cohomology class $\tau_{[N,E,\Phi,c]}$ is
\beq
\lefteqn{\langle \tau_{[N,E,\Phi,c]}\, , \, \Ind([P]) \rangle  =}\non\\
&& \langle [N\times_B S^*_\pi M \, , \, E\times_B S^*_\pi M \, , \, \pi_*^{-1}(\Phi) \, , \, \Td(T_\pi M\otimes\cc)\wedge \pi^*(c)] \, , \, [P] \rangle \non
\eeq
where $[P]\in K_1(\cinfc(S^*_\pi M)\rtimes G)$ is the leading symbol class of $P$. \hfill\cqfd
\end{corollary}


\section{Foliated dynamical systems}\label{sfol}

Let $(V,\Fc)$ be a compact foliated manifold without boundary. We recall \cite{C94} that its \emph{holonomy groupoid} $H$ is a Lie groupoid with $V$ as set of units, and the arrows are equivalence classes of leafwise paths $\gamma$ from a source point $x=s(\gamma)$ to a range point $y=r(\gamma)$ belonging to the same leaf, with the following relation: two leafwise paths $\gamma$ and $\gamma'$ from $x$ to $y$ are equivalent if and only if they induce the same diffeomorphism (holonomy) from a transverse neighborhood of $x$ to a transverse neighborhood of $y$. The composition of arrows is the concatenation product of paths. It is well-known that the holonomy groupoid can be reduced to a Morita equivalent \emph{\'etale} groupoid upon a choice of complete transversal. Indeed, let $B\to V$ be an immersion of a closed manifold everywhere transverse to the leaves and intersecting each leaf at least once. The subgroupoid 
\be
H_B=\{ \gamma\in H \ |\ s(\gamma)\in B\ \mbox{and}\ r(\gamma)\in B\}
\ee
is \'etale, i.e. the range and source maps from $H_B$ to $B$ are local diffeomorphisms, and is Morita equivalent to the holonomy groupoid. $H_B$ naturally acts on the submersion $s:M\to B$, corresponding to the restriction of the source map $s:H\to V$ to the submanifold
\be
M=\{ x\in H\ |\ s(x)\in B\}\ .
\ee
The right action of $H_B$ on $M$ is given by composition of arrows in $H$: for any $x\in M$ and $\gamma\in H_B$ such that $s(x)=r(\gamma)\in B$, the composite $x\cdot \gamma$ is indeed in $M$. The range map $r:M\to V$, $x\mapsto r(x)$ is a covering, mapping the fibers of the submersion $M$ to the leaves of $V$.\\ 
Now suppose in addition the foliation $(V,\Fc)$ endowed with a transverse flow of $\rr$. Hence there is a one-parameter group of diffeomorphisms $\phi_t$, $t\in \rr$ on $V$, such that $\phi_t$ maps leaves to leaves, and the vector field $\Phi$ generating the flow is nowhere tangent to the leaves (thus in particular does not vanish). Since $M\to V$ is a covering, the generator of the flow $\phi$ can be lifted in a unique way to a vector field also denoted $\Phi$ on the submersion $M$. Since $\phi$ preserves the leaves, $\Phi$ descends to a vector field $\Phib$ on the base $B$. By hypothesis $B$ is a closed manifold, hence $\Phib$  can be integrated to a flow $\phib$ on $B$, and consequently $\Phi$ also generates a flow $\phi$ on $M$, mapping fibers to fibers. By construction the range map $r:M\to V$ is $\rr$-equivariant:
\be
r(\phi_t(x)) = \phi_t(r(x))\quad \forall\  x\in M\ ,\ t\in \rr\ .
\ee
The vector field $\Phib$ on $B$ is obtained as follows: at each point $b$ of the submanifold $B\subset V$, the tangent space $T_bV$ decomposes canonically as the direct sum of $T_bB$ and the tangent space to the leaf. $\Phi$ projects accordingly to a vector field on $B$ precisely corresponding to the generator $\Phib$. The flow $\phib$ also lifts to a flow on the groupoid $H_B$ because the range and source maps $H_B\rightrightarrows B$ are \'etale. This results in an action of $\rr$ on the groupoid $H_B$ by homomorphisms, i.e. $\phib_t(\gamma\delta)=\phib_t(\gamma)\phib_t(\delta)$ for all $\gamma,\delta\in H_B$ and $t\in \rr$. We form a new Lie groupoid by taking the crossed product
\be
G = H_B\cp\rr\ .
\ee
The arrows of $G$ are pairs $(\gamma,t)\in H_B\times \rr$, with composition law $(\gamma,t)(\delta,u) = (\gamma\, \phib_{-t}(\delta),t+u)$. Combining the respective actions of $H_B$ and $\rr$ on the submersion $s:M\to B$ yields the following action of the crossed product: $x\cdot (\gamma,t) = \phi_t(x\cdot\gamma) = \phi_t(x)\cdot \phib_t(\gamma)$. To summarize we are left with two Morita equivalent groupoids, namely $V\cp\rr$ and $M\cp G$. This equivalence can be realized at the level of convolution algebras via a homomorphism 
\be
\rho:\cinfc(V\cp\rr)\to \cinfc(M\cp G)
\ee 
constructed as follows. The range map $r:M\to V$ is a covering of a compact manifold, hence choosing and partition of unity relative to a local trivialization of the covering one can build a ``cut-off'' function $c\in \cinfc(M)$ with the property $\sum_{x\in r^{-1}(v)}c(x)^2 = 1$ for all $v\in V$. Then $\rho$ sends $f\in \cinfc(V\cp \rr)$ to the function $\rho(f)\in \cinfc(M\cp G)$ defined by
\be
\rho(f)(x,(\gamma,t)) = c(x)\, f(r(x), t)\, c(\phi_t(x\cdot\gamma))
\ee
for any $x\in M$, $\gamma\in H_B$ and $t\in \rr$. Note that $\phi_t(x\cdot\gamma)$ is the action of $(\gamma,t)\in G$ on $x$. One checks that $\rho$ is a homomorphism of algebras, and that the resulting map in $K$-theory $\rho_!: K_0(\cinfc(V\cp\rr)) \to K_0(\cinfc(M\cp G))$ is independent of any choice of cut-off function.\\ 

Now let $E=E_+\oplus E_-$ be a $\zz_2$-graded, $\rr$-equivariant vector bundle over $V$. We consider an $\rr$-invariant \emph{leafwise} differential elliptic operator $D$ of order one and odd degree acting on the sections of $E$. Hence according to the $\zz_2$-grading $D$ splits as the sum of two differential operators of order one acting only in the directions of the leaves $D_+:\cinfc(V,E_+)\to\cinfc(V,E_-)$ and $D_-:\cinfc(V,E_-)\to\cinfc(V,E_+)$. This may be rewritten in the usual matrix form
\be
D= \left( \begin{matrix} 
0 & D_- \\
D_+ & 0 \end{matrix} \right)\ .
\ee
We make the further assumption that $D$ is formally self-adjoint with respect to some metric on $V$ and hermitean structure on $E$. Note that the condition of $\rr$-invariance for $D$ imposes strong restrictions on the flow on $V$. The vector bundle $E$ defines a vector bundle (still denoted by $E$) on the manifold $M$ by pullback  with respect to the range map $r$. Since $D$ is an $\rr$-invariant leafwise differential operator on $V$, it pullbacks accordingly to a $G$-invariant fiberwise differential operator on the submersion $s:M\to B$. Let $S^1=\rr/\zz$ be the circle, endowed with the trivial action of $G$. We view the product $M'=S^1\times M$ as a submersion over the base manifold $B$, endowed with its $G$-action, and consider the differential operator
\be
Q= \left( \begin{matrix} 
\d_x & D_- \\
D_+ & -\d_x \end{matrix} \right)
\ee
where $\d_x=\frac{\d}{\d x}$ denotes partial differentiation with respect to the circle variable $x\in [0,1]$. Then $Q$ is a $G$-invariant fiberwise differential operator on $M'$, and  is elliptic because $Q^2=D^2 + (\d_x)^2$. As pseudodifferential operators, the modulus $|Q|$ and its parametrix $|Q|^{-1}$ are well-defined only modulo the ideal of smoothing operators. Hence the ``sign'' of $Q$ and the ``spectral projection'' onto the 1-eigenspace 
\be
F=\sign(Q)=\frac{Q}{|Q|}\ ,\qquad P=\frac{1+F}{2}
\ee 
are represented by fiberwise pseudodifferential operators of order zero on the submersion $M'\to B$, and modulo smoothing operators $F$ and $P$ are $G$-invariant and fulfill the identities $F^2=1$ and $P^2=P$. Let $[e]\in K_0(\cinfc(V\cp\rr))$ be a $K$-theory class represented by an idempotent matrix $e$. For simplicity we assume $e\in \cinfc(V\cp \rr)$. The suspension of the idempotent $\rho(e)\in \cinfc(M\cp G)$ is the invertible element
\be
u=1 + \rho(e)(\beta-1)\ \in \cinfc(M'\cp G)^+\ ,
\ee
where the function $\beta\in \cinf(S^1)$, $\beta(x)=\exp(2\pi i\, x)$  is the Bott generator of the circle. The algebra $\cinfc(M')$ acting by pointwise multiplication on the space of sections $\cinfc(M',E)$ is naturally represented in the algebra of sections of the bundle of pseudodifferential operators of order zero $\CL^0_c(M',E)$. This representation is $G$-equivariant, hence extends to a representation of $\cinfc(M')\cp G$ in the convolution algebra $\Ec = \cinfc(B,\CL^0_c(M',E))\cp G$. Accordingly $u$ is represented by an invertible element  $U$ in the unitalization of $\Ec$. The Toeplitz operator
\be
T= PUP + (1-P)\ \in \Ec^+
\ee
is uniquely defined modulo smoothing operators. One has $T\equiv 1+P(U-1)$ modulo the ideal $\Bc \subset \Ec$ of order $-1$ pseudodifferential operators, and the inverse of $T$ modulo this ideal is represented by the Toeplitz operator $PU^{-1}P + (1-P)\equiv 1+ P(U^{-1}-1)$. The quotient $\Ec/\Bc$ is isomorphic to the convolution algebra of leading symbols $\Ac = \cinfc(B,\LS^0_c(M',E))\cp G$. Note that the leading symbol of $T$ is uniquely defined and reads
\be
\si_T = 1 + \si_P (u-1)\ \in \Ac^+\ ,
\ee
where the leading symbol $\si_P$ of $P$ is a $G$-invariant idempotent in the algebra of sections of $\LS^0(M',E)$.  Hence the equality $\si_Pu=u\si_P$ holds in $\Ac$ and the inverse of $\si_T$ is exactly $1+\si_P(u^{-1}-1)$. It is easy to see that the $K$-theory class $[\si_T]\in K_1(\Ac)$ only depends on the $K$-theory class $[e]\in K_0(\cinfc(V\cp\rr))$ of the idempotent $e$. 

\begin{definition}\label{dind}
Let $(V,\Fc)$ be a compact foliated manifold endowed with a transverse flow of $\rr$, and $D$ be an $\rr$-equivariant leafwise elliptic differential operator of order one and odd degree acting on the sections of a $\zz_2$-graded equivariant vector bundle over $V$. Then for any class $[e]\in K_0(\cinfc(V\cp\rr))$, we define the index of $D$ with coefficients in $[e]$ as the $K$-theory class
\be
\Ind(D,[e]) = \Ind_E([\si_T]) \in  K_0(\Bc)\ ,
\ee
where $T$ is the Toeplitz operator constructed above and $\Ind_E$ is the index map of the extension $(E)\ :\ 0\to \Bc \to \Ec \to \Ac \to 0$.
\end{definition}

From now on we specialize to a \emph{codimension one} foliation $(V,\Fc)$ endowed with a transverse flow $\phi$. The connected components of the closed transversal $B$ are all diffeomorphic to the circle. Hence, the induced flow $\phib$ on a connected component $B_p$ acts by rotation with a given period $p\in \rr$. Choosing a base-point we canonically get a parametrization of $B_p$ by the variable $b\in \rr/p\zz$, and the flow is simply $\phib_t(b) \equiv b+t$ mod $p$. In the same way, we consider the range $r(\gamma)\in B_p$ and the source $s(\gamma)\in B_q$ of any arrow $\gamma\in H_B$ as numbers respectively in $\rr/p\zz$ and $\rr/q\zz$, and denote by $r^*db$, $s^*db$ the measures on $H_B$ pulled back from the Lebesgue measure by the \'etale maps $r$, $s$ respectively. Note that $db$ provides a canonical holonomy-invariant transverse measure on $(V,\Fc)$ because the flow $\phi$ preserves the leaves of the foliation. In particular $r^*db=s^*db$ on $H_B$. Remark also that because the transverse Lebesgue measure is holonomy-invariant, the holonomy of a closed leafwise path is always the identity. As a consequence, the holonomy groupoid $H$ and its reduced groupoid $H_B$ have no automorphisms except the units (that is, $r(\gamma)=s(\gamma)$ implies $\gamma$ is a unit).

\begin{proposition}\label{ptr}
The convolution algebra $\cinfc(G)$ acts on the space of functions $\cinf(B)$ by smoothing operators. In this representation, the operator trace of an element $f\in \cinfc(G)$ reads
\be
\Tr(f)= \sum_{B_p} \sum_{n\in\zz} \int_{H_B} f(\gamma,np+ r(\gamma)-s(\gamma))\, r^*db(\gamma)\label{op}
\ee
where the sum runs over all connected components $B_p$ of the flow $\phib$ on the transversal, $p$ is the period of $B_p$ and $db$ is the pullback of the Lebesgue measure on the groupoid $H_{B_p}$. In particular the distribution kernel of $\Tr$ is a measure on $G$ whose support is the manifold of automorphisms in the groupoid $G$.
\end{proposition}
{\it Proof:} We define the action of $f\in \cinfc(G)$ on a test function $\xi\in \cinf(B)$ as follows: at any point $b\in B$,
$$
(f\cdot\xi)(b) = \sum_{\gamma\in H_B^b}\int_{-\infty}^{\infty} f(\gamma,t) \xi(\phib_t(s(\gamma)))\, dt
$$
where $H_B^b$ is the set of arrows $\gamma\in H_B$ such that $r(\gamma)=b$. The sum is finite and the integral converges because $f$ is of compact support on $G$. One readily verifies that this defines a representation of the convolution algebra $\cinfc(G)$ on $\cinf(B)$, and that the distributional kernel of the operator $f$ on the manifold $B\times B$ is 
$$
k_f(b,b') = \sum_{\gamma\in H_B^b} \int_{-\infty}^{+\infty} f(\gamma,t) \, \delta(\phib_t(s(\gamma)) - b')\, dt
$$
where $\delta$ is the Dirac measure. If $s(\gamma)$ is in a connected component of period $p$ for the flow $\phib$, one has $\phib_t(s(\gamma))\equiv s(\gamma)+t$ mod $p$ and the integral selects all values $t=np+b'-s(\gamma)$, $n\in \zz$. Therefore
$$
k_f(b,b') = \sum_{n\in\zz}\, \sum_{\gamma\in H_B^b} f(\gamma, np +b' -s(\gamma))
$$
(where $p$ is the period of the connected component of $s(\gamma)$) is a smooth function. Hence $\cinfc(G)$ acts by smoothing operators. The trace $\Tr(f)$ is given by the integral, with respect to the measure $db$, of the kernel restricted to the diagonal $k_f(b,b)$. In this case $b=b'=r(\gamma)$ and $s(\gamma)$ belong to the same connected component, say $B_p$. This leads to
\beq
\Tr(f) &=& \sum_{B_p} \sum_{n\in\zz} \int_0^p \sum_{\gamma\in H_{B_p}^b} f(\gamma, np +b -s(\gamma)) \, db \non\\
&=& \sum_{B_p} \sum_{n\in\zz} \int_{H_{B_p}} f(\gamma,np+ r(\gamma)-s(\gamma))\, r^*db(\gamma) \non
\eeq
as claimed. If $b=r(\gamma)\in B_p$ and $t=np+ r(\gamma)-s(\gamma)$ for some $n\in \zz$, one has $\phib_t(s(\gamma)) \equiv r(\gamma)$ mod $p$ which is equivalent to $b\cdot(\gamma,t)= b$, that is $(\gamma,t)\in G$ is an automorphism.  \cqfd\\

The operator trace on $\cinfc(G)$ is therefore a cyclic zero-cocycle localized at the isotropic set in $G$, and is of order $k=0$ because its distribution kernel is a measure. We may split it into several parts. Let us first consider the linear functional $\Tr_0:\cinfc(G)\to\cc$ whose support is localized at the set of \emph{units} $B\subset G$. It amounts to sum only over the arrows $\gamma=r(\gamma)=s(\gamma)\in B$ and the integer $n=0$ in the right-hand-side of (\ref{op}):
\be
\Tr_0(f) = \int_B f(b,0) \, db\ .
\ee
One easily checks that $\Tr_0$ is a trace. By construction it is localized at the submanifold of units in $G$, and its distribution kernel corresponds to the (holonomy invariant) transverse measure $db$. There is an associated cyclic cocycle $\tau_0$ over the algebra $\Bc$.

\begin{definition}
Let $\Tr_0:\cinfc(G)\to\cc$ be the trace localized at units given by the holonomy invariant transverse measure on $(V,\Fc)$ induced by the flow $\phi$, and $\tau_0$ be the associated cyclic cocycle over $\Bc$. We define the \emph{Connes-Euler characteristics} of the $\rr$-invariant operator $D$ with coefficients in a $K$-theory class $[e]\in K_0(\cinfc(V\cp\rr))$ as the complex number
\be
\chi(D,[e]) = \langle \tau_0, \Ind(D,[e]) \rangle 
\ee
\end{definition}

This number is the analogue, in the $\rr$-equivariant context, of the index defined by Connes for longitudinal elliptic operators on a foliation endowed with a holonomy invariant transverse measure \cite{C94}. Now let $\tau$ be the cyclic cocycle over $\Bc$ corresponding to the full operator trace $\Tr$. Our aim is to evaluate $\tau$ on the $K$-theory class $\Ind(D,[e])\in K_0(\Bc)$. As we shall see, the complementary part $\langle \tau, \Ind(D,[e]) \rangle - \chi(D,[e])$ is related to the \emph{periodic orbits} of the flow $\phi$ on $(V,\Fc)$. If $\Pi\subset V$ is a periodic orbit with period $p_{\Pi}$, a choice of base-point canonically provides a parametrization of the orbit by a variable $v\in [0,p_{\Pi}]$, with $\phi(v)\equiv v+t$ mod $p_{\Pi}$. The flow at $t=p_{\Pi}$ is called the \emph{return map}. It defines an endomorphism $h_{\Pi}'(v)$ on the tangent space to the leaf at any point $v\in \Pi$, together with an even-degree endomorphism $j_{\Pi}(v)$ on the fiber of the $\zz_2$-graded vector bundle $E$ over $v$. We shall suppose that $h_{\Pi}'$ is non-degenerate, in the sense that $\det (1-h_{\Pi}'(v))\neq 0$ for all $v\in \Pi$. This implies in particular that the periodic orbits of the flow are isolated.

\begin{proposition}\label{ptheta}
Let $\Pi\subset V$ be a periodic orbit of the flow $\phi$, and $v\in \Pi$. We denote by $h_{\Pi}'(v)$ the action of the return map on the leafwise tangent space, and by $\tr_s(j_{\Pi}(v))$ the supertrace of the return map on the $\zz_2$-graded vector bundle $E$ at $v$. Suppose that $h_{\Pi}'$ is non-degenerate. Then the linear functional $\Theta_{\Pi}:\cinfc(V\cp\rr)\to\cc$ 
\be
\Theta_{\Pi}(f) = \sum_{n\in \zz^*} \int_{\Pi} \frac{\tr_s(j_{\Pi}(v)^n)}{|\det(1-h_{\Pi}'(v)^n)|}\,f(v,np_{\Pi})\, dv  \label{theta}
\ee
with $p_{\Pi}\in \rr$ the period of the orbit, is a trace. Actually the function $v\mapsto \tr_s(j_{\Pi}(v)^n)/|\det(1-h_{\Pi}'(v)^n)|$ is constant along the orbit.
\end{proposition}
{\it Proof:} If $v$ and $v'$ are two distinguished points in the orbit $\Pi$, the endomorphisms $j_{\Pi}(v)$ and $j_{\Pi}(v')$ are conjugate, hence have the same supertrace. Similarly $\det(1-h'_{\Pi}(v))$ and $\det(1-h'_{\Pi}(v'))$ are equal. The function $v\mapsto \tr_s(j_{\Pi}(v)^n)/|\det(1-h_{\Pi}'(v)^n)|$ is therefore constant along the orbit, and can be pushed out of the integration over $v$. In fact for any fixed $n\in\zz^*$, the integral
$$
\Te_n(f) = \int_{\Pi}f(v,np_{\Pi})\, dv
$$
defines a trace on $\cinfc(V\cp\rr)$. Indeed for two functions $f,g$ on $\Pi\times \rr$ one has 
\beq
\Theta_n(fg) &=& \int_{\Pi} \int_{-\infty}^{\infty} f(v,np_{\Pi}-t)g(v-t,t)dt\, dv \non\\
&=& \int_{\Pi} \int_{-\infty}^{\infty} f(v,t)g(v+t,np_{\Pi}-t)dt\, dv \non\\
&=& \int_{\Pi} \int_{-\infty}^{\infty} f(v-t,t)g(v,np_{\Pi}-t)dt\, dv\ =\ \Theta_n(gf)\non
\eeq
where we used repeatedly the fact that the functions $f$ and $g$ are $p_{\Pi}$-periodic in the variable $v$. \cqfd\\

\begin{theorem}\label{tfolp}
Let $(V,\Fc)$ be a codimension one compact foliated manifold endowed with a transverse flow, $E$ a $\zz_2$-graded $\rr$-equivariant vector bundle over $V$, and $D$ an odd $\rr$-invariant leafwise elliptic differential operator of order one. Assume that the periodic orbits of the flow are non-degenerate. Then for any class $[e]\in K_0(\cinfc(V\cp \rr))$, the pairing of the index $\Ind(D,[e])$ with the cyclic cocycle $\tau$ induced by the operator trace is
\be
\langle [\tau] , \Ind(D,[e]) \rangle =\chi(D,[e]) +  \sum_{\Pi}\Theta_{\Pi}(e)\ ,
\ee
where $\chi(D,[e])$ is the Connes-Euler characteristics, the sum runs over the periodic orbits $\Pi$ of the flow, and $\Theta_{\Pi}$ is the canonical trace (\ref{theta}) on $\cinfc(V\cp \rr)$.  
\end{theorem}
{\it Proof:} We use zeta-function renormalization by means of the complex powers of the operator $|Q|$. Since $\Tr$ is a trace on $\cinfc(G)$ of order zero and localized at the submanifold of automorphisms, Eq. (\ref{ctr}) implies
$$
\langle [\tau] , \Ind(D,[e]) \rangle = C\circ\Res( T^{-1} [\ln|Q|,T])
$$
with $T^{-1}$ the inverse of $T$ modulo smoothing operators, and $C$ the distribution kernel of $\Tr$ which is a measure on $G$. Note that $T^{-1}[\ln|Q|,T]$ is a pseudodifferential operator of order $-1$ because $|Q|$ is $G$-invariant modulo smoothing operators. We know that $C$ splits as a sum of a part $C_0$ whose support is localized at units, and a complementary part $C_p$. By definition of the Connes-Euler characteristics,
$$
\langle [\tau] , \Ind(D,[e]) \rangle - \chi(D,[e]) = {C_p}\circ\Res( T^{-1} [\ln|Q|,T])\ .
$$
The support of $C_{\pi}$ is the set of arrows $(\gamma,t)\in G$ which are not units and such that $r(\gamma)=\phib_t(s(\gamma))$ in $B$. Such $(\gamma,t)$ acts on the fiber $M_{r(\gamma)}$ of the submersion $s:M\to B$ by a diffeomorphism $h:x\mapsto \phi_t(x\cdot\gamma)$. Since $(\gamma,t)$ is not a unit, one sees that $x\in M_{r(\gamma)}$ is a fixed point of $h$ if and only if its image $r(x)\in V$ is contained in a periodic orbit of the flow. By hypothesis the periodic orbits are isolated, so the fixed points for $h$ is discrete subset of $M_{r(\gamma)}$. Extending the action of $h$ on the fiber $M'_{r(\gamma)} = S^1\times M_{r(\gamma)}$ in a trivial way, the fixed point set becomes a discrete union of circles. In this simple situation, formula (\ref{floc}) reduces to 
$$
\res \Tr(RU_h|Q|^{-z}) = \sum_{\substack{S^1\subset M'_{r(\gamma)} \\ \mathrm{fixed}}}\int_{S^*S^1} \frac{[\si_R]_{-1}(x,\xi;0,0)}{|\det(1-h')|} \, \frac{\xi dx}{2\pi}
$$
for any pseudodifferential operator $R$ of order $-1$ on $M'_{r(\gamma)}$. Here $(x,\xi)$ are the canonical coordinates on the cotangent bundle of $S^1$, $[\si_R]_{-1}$ is the leading symbol of $R$, and $h'$ is the differential of the diffeomorphism $h$ in the direction normal to the circle, i.e. $h'$ is an endomorphism of the tangent space to $M_{r(\gamma)}$ at the fixed points. We apply this formula to the operator $R=T^{-1} [\ln|Q|,T](\gamma,t)$. At the level of leading symbols, one has
$$
\si_{T^{-1} [\ln|Q|,T]}(x,\xi;0,0) \sim \si_T^{-1} \frac{\d \si_T}{\d x}(x,\xi;0,0)
$$
where $\si_T = 1+\si_P(u-1)$ is the leading symbol of $T$, $u=1+\rho(e)(\beta-1)$, and $\si_P$ is the leading symbol of the spectral projection $P$. In matrix notation adapted to the $\zz_2$-grading of the vector bundle $E$,
$$
\si_P(x,\xi;0,0) = \left( \begin{matrix} \te(\xi) & 0 \\ 
                                  0  & \te(-\xi) \end{matrix} \right)\ ,
$$
where $\te: \rr\to\{0,1\}$ is the Heavyside function. The dependency of $\si_T$ upon the parameter $x$ comes only from the Bott element $\beta\in \cinf(S^1)$. The cosphere bundle $S^*S^1$ being just a sum of two copies of the circle, the computation of the residue is straightforward. It is a sum over the (isolated) fixed points of $h$ in $M_{r(\gamma)}$:
$$
\res \Tr(T^{-1} [\ln|Q|,T](\gamma,t)U_h|Q|^{-z}) = \sum_{\substack{y\in M_{r(\gamma)} \\ \mathrm{fixed}}} \frac{\tr_s(\rho(e)(y,(\gamma,t)))}{|\det(1-h'(y))|}
$$
Since $y\in M$ is fixed by $(\gamma,t)\in G$, one has $\rho(e)(y,(\gamma,t)) = c(y)^2 e(r(y),t)$ where $c\in\cinfc(M)$ is the cut-off function used in the construction of the homomorphism $\rho$, and $r(y)\in V$ is the projection of $y$. Since the \'etale groupoid $H_B$ acts without fixed points on $B$, integrating with respect to the measure $C_p$ on $G$ amounts to integrate over all the periodic orbits $\Pi\subset V$ of the flow. A straightforward computation yields
$$
{C_p}\circ\Res( T^{-1} [\ln|Q|,T]) = \sum_{\Pi} \sum_{n\in\zz^*}\int_{\Pi} \frac{\tr_s(j_{\Pi}(v)^n)}{|\det(1-h_{\Pi}'(v)^n)|}\,e(v,np_{\Pi})\, dv
$$
where $h'_{\Pi}(v)$ and $j_{\Pi}(v)$ are the actions of the return map respectively on the tangent space to the leaf and on the fiber of $E$ at $v\in \Pi$. \cqfd

\begin{corollary}
If the Connes-Euler characteristics $\chi(D,[e])$ is not an integer, then the flow $\phi$ has periodic orbits.
\end{corollary}
{\it Proof:} Indeed the cyclic cocycle $\tau$ over $\Bc$ comes from the operator trace on $\cinfc(G)$, so the pairing $\langle [\tau] , \Ind(D,[e]) \rangle$ is always an integer. The sum over periodic orbits cannot vanish if the Connes-Euler characteristic does not belong to $\zz$. \cqfd\\

We end this section with a brief discussion of the more general situation of a flow $\phi$ compatible with the foliation $(V,\Fc)$, i.e. mapping leaves to leaves, but also admitting fixed points. Remark that a leaf containing a fixed point is necessarily preserved by $\phi_t$ for all $t\in \rr$: in this case the vector field $\Phi$ generating the flow is tangent to the leaf. We shall suppose that the flow $\phi$ is \emph{non-degenerate} in the folIowing sense:\\

\noindent i) At any fixed point $v\in V$, the tangent map $T_v\phi_t$ on the tangent space $T_vV$ has no singular value equal to $1$ for $t\neq 0$;\\

\noindent ii) Each periodic orbit $\Pi\subset V$ is transverse to the leaves, and the return map induced by the flow on the tangent space to the leaf $T_vL$ at any $v\in \Pi$ has no singular value equal to $1$.\\

\noindent iii) If $v\in V$ belongs to a leaf $L$ preserved by the flow, the tangent map $T_v\phi_t$ acting on the transverse space $T_vV/T_vL$ has no singular value equal to $1$ for $t\neq 0$;\\

\noindent iv) The holonomy of a closed path contained in a leaf $L$ preserved by the flow is always the identity.\\

\noindent Conditions i) and ii) immediately imply that the fixed points and the periodic orbits of the flow are isolated. In particular, a fixed point cannot belong to a dense leaf. If $B\subset V$ is a complete closed transversal, the induced flow $\phib$ on $B$ may also have fixed points, consisting in the intersection of $B$ with all the leaves preserved by $\phi$. Conditions i) and ii) also imply that the fixed points of $\phib$ are isolated in $B$. Moreover,  by condition iii) the tangent map $T_b\phib_t$ acting on the tangent space $T_bB$ at a fixed point $b\in B$ is not the identity if $t\neq 0$; thus $T_b\phib_t$ acts by the multiplication operator $e^{\kappa_bt}$ where $\kappa_b\neq 0$ is the \emph{exponent} of the flow at $b$. We add condition iv) because, unlike in the previous situation of a transverse flow, the holonomy of a closed path is not automatically the identity if the leaf is preserved by the flow. Now the closed manifold $B$ can be partitioned into the orbits of the flow $\phib$, which are of three different types:
\begin{itemize}
\item The periodic orbits, denoted $B_p$, all diffeomorphic to a circle. Any such orbit of period $p$ is naturally parametrized by $b\in \rr/p\zz$ once a base-point is chosen, and the flow is $\phib_t(b)\equiv b+t$ mod $p$ for all $t$. We denote $db$ the Lebesgue measure on $B_p$;
\item The infinite orbits, denoted $B_{\infty}$, all diffeomorphic to $\rr$. They are also naturally parametrized by $b\in\rr$ once a base-point is chosen, with $\phib_t(b)=b+t$ for all $t$. We denote $db$ the Lebesgue measure on $B_{\infty}$;
\item The fixed points $\phib_t(b)=b$ for all $t$, each of them connecting two infinite orbits. 
\end{itemize}
The reduced holonomy groupoid $H_B$ still carries an action of $\rr$ by homomorphisms, so we can form the crossed product $H_B\cp\rr$ as before. If $B_p$ is a periodic orbit, we denote $H_{B_p}$ the subgroupoid of $H_B$ consisting in arrows $\gamma$ with source and range contained in $B_p$. Similarly with the infinite orbits. As before, the measure $db$ on periodic and infinite orbits is holonomy-invariant. However, compared with the previous situation of a transverse flow, the main difficulty with the presence of fixed points is that the convolution algebra $\cinfc(G)$ does no longer act by smoothing operators on the space $\cinf(B)$. In order to define a trace we are forced to restrict to the subalgebra
\be
\cinfc(G_+^*)=\{f\in\cinfc(G)\ |\ \supp(f)\subset H_B\times \rr_+^*\}
\ee
where $\rr_+^*\subset \rr$ is the abelian monoid (with addition) of real numbers $>0$. One easily checks that the convolution product on $\cinfc(G)$ preserves the subspace $\cinfc(G_+^*)$. The submanifold $G_+^*=H_B\times\rr_+^*$ of $G$ has a partially defined composition law coming from the composition of arrows in $G$; however $G_+^*$ is not a subgroupoid. We nevertheless regard $\cinfc(G_+^*)$ as the convolution algebra of the ``crossed-product" $G_+^*=H_B\cp\rr_+^*$. 

\begin{proposition}
The algebra $\cinfc(G_+^*)$ acts on the space $\cinf(B)$ by distributional kernels which may be singular on the diagonal. When the flow $\phi$ is non-degenerate, the operator trace on smoothing kernels canonically extends to a trace on $\cinfc(G_+^*)$ by
\beq
\Tr(f) &=& \sum_{B_p} \sum_{n\in\zz} \int_{H_{B_p}} f(\gamma,np+r(\gamma)-s(\gamma)) \, r^*db(\gamma) \non\\
&& + \sum_{B_{\infty}} \int_{H_{B_{\infty}}} f(\gamma,r(\gamma)-s(\gamma)) \, r^*db(\gamma) \label{trace} \\
&& + \sum_b  \int_0^{\infty} \frac{f(b,t)}{|1-e^{\kappa_bt}|}\, dt \non
\eeq
where the sums run over all periodic orbits $B_p$, infinite orbits $B_{\infty}$ and fixed points $b$ of the flow $\phib$ on the transversal; $p>0$ is the period of the orbit $B_p$ and $\kappa_b\neq 0$ is the exponent of $\phib$ at the fixed point $b$.
\end{proposition}
{\it Proof:} As in the proof of Proposition \ref{ptr}, the distributional kernel of an element $f\in\cinfc(G_+^*)$ acting on $\cinf(B)$ is 
$$
k_f(b,b') = \sum_{\gamma\in H_B^b} \int_0^{\infty} f(\gamma,t) \, \delta(\phib_t(s(\gamma)) - b')\, dt\ .
$$
Its restriction to the diagonal is still a smooth function except at the fixed points of the flow $\phib$ where it is proportional to a Dirac measure. Its integral over $B$ thus extends the operator trace of smoothing kernels. We split $B$ into three parts corresponding to the different types of connected components of the flow. The first term of (\ref{trace}) associated to the periodic components is calculated exactly as in \ref{ptr}, and similarly for the second term associated to the infinite components. At a fixed point, the equality $\phib_t(s(\gamma))=r(\gamma)$ implies $s(\gamma)=r(\gamma)$, hence $\gamma$ is necessarily a unit by hypothesis iv). The kernel $k_f(b,b)$ for $b$ in the vicinity of a fixed point thus reads
$$
k_f(b,b) = \int_0^{\infty} f(b,t)\, \delta(\phib_t(b)-b)\, dt \ .
$$
Taking the fixed point as the origin of a local coordinate system, the equality $\phib_t(b)\sim b e^{\kappa t}$ holds at first order in the variable $b$, where $\kappa$ is the exponent of $\phib$ at the fixed point. Thus
$$
k_f(b,b)= \delta(b) \int_0^{\infty} \frac{f(b,t)}{|1-e^{\kappa t}|}\, dt \ ,
$$
and integrating over $b$ yields the contribution of the fixed point to the trace, whence the third term of (\ref{trace}).  \cqfd\\

Let $\cinfc(V\cp\rr_+^*)$ be the space of smooth compactly supported functions on $V\times\rr$ with support contained in the open subset $V\times \rr_+^*$. It is a subalgebra of the convolution algebra $\cinfc(V\cp\rr)$. Let $D$ be an $\rr$-equivariant leafwise elliptic differential operator $D$ of order one and odd degree acting on the sections of a $\zz_2$-graded equivariant vector bundle $E$ over $V$, and $e\in\cinfc(V\cp\rr_+^*)$ be an idempotent. The index $\Ind(D,[e])$ defined in \ref{dind} is then a $K$-theory class of the subalgebra $\Bc_+^*=\cinfc(B,\CL_c^{-1}(M',E))\cp G_+^*$ of $\Bc$. As before we denote $\tau$ the cyclic cocycle over $\Bc_+^*$ induced by the trace $\Tr:\cinfc(G_+^*)\to\cc$. The pairing $\langle [\tau] , \Ind(D,[e]) \rangle$ is a well-defined complex number. Its explicit computation involves, besides the periodic orbits of the flow, also the fixed points.

\begin{proposition}\label{pW}
Let $v\in V$ be a fixed point of the non-degenerate flow $\phi$. We denote by $\kappa_v$ the generator of the flow on the tangent space $T_vV$, and by $j_v$ the generator of the flow on the fiber $E_v$. Then the linear functional $W_v:\cinfc(V\cp\rr_+^*)\to\cc$ 
\be
W_v(f) =  \int_0^{\infty} \frac{\tr_s(e^{j_vt})}{|\det(1-e^{\kappa_vt})|}\,f(v,t)\, dt \label{W}
\ee
is a trace. 
\end{proposition}
{\it Proof:} Since $\phi_t(v)=v$, $\forall t>0$, one has for any functions $f,g\in \cinfc(V\cp\rr_+^*)$
$$
(fg)(v,t) = \int_{-\infty}^{+\infty} f(v,u)g(v,t-u)\, du = (gf)(v,t)
$$
so that $W_v$ is obviously a trace. \cqfd\\

Since the support of the index $\Ind(D,[e]) \in K_0(\Bc_+^*)$ does not meet the part $t\leq 0$ of the groupoid $G$, the Connes-Euler characteristics does not appear in the computation of $\langle [\tau] , \Ind(D,[e]) \rangle$. One is left with the periodic orbits and the fixed points only. The method of proof of Theorem \ref{tfolp} leads at once to the following result.

\begin{theorem}
Let $(V,\Fc)$ be a codimension one compact foliated manifold endowed with an $\Fc$-compatible flow, $E$ a $\zz_2$-graded $\rr$-equivariant vector bundle over $V$, and $D$ an odd $\rr$-invariant leafwise elliptic differential operator of order one. Assume that the periodic orbits and the fixed points of the flow are non-degenerate. Then for any class $[e]\in K_0(\cinfc(V\cp \rr_+^*))$, the pairing of the index $\Ind(D,[e])$ with the cyclic cocycle $\tau$ induced by the operator trace is
\be
\langle [\tau] , \Ind(D,[e]) \rangle = \sum_{\Pi}\Theta_{\Pi}(e) + \sum_{v\ \mathrm{fixed}} W_v(e)\ ,
\ee
where the sums runs over the periodic orbits $\Pi$ and the fixed points $v$ of the flow, $\Theta_{\Pi}$ and $W_v$ are the canonical traces (\ref{theta}) and (\ref{W}) on $\cinfc(V\cp \rr_+^*)$. 
\end{theorem}
{\it Proof:} As in the proof of \ref{tfolp} we let $C$ be the measure on $G_+^*$ correponding to the trace $\Tr:\cinfc(G_+^*)\to\cc$. Then  
$$
\langle [\tau] , \Ind(D,[e]) \rangle = C\circ\Res( T^{-1} [\ln|Q|,T])
$$
According to (\ref{trace}), the trace can be decomposed into three parts respectively associated to periodic orbits, infinite orbits and fixed points of the flow $\phib$ on the transversal $B$. Hence $C$ can be accordingly decomposed into three measures $C_p$, $C_{\infty}$ and $C_f$ over $G_+^*$. The contributions of $C_p$ and $C_{\infty}$ to the pairing $\langle [\tau] , \Ind(D,[e]) \rangle$ are computed exactly as in \ref{tfolp} and yield the sum over the periodic orbits $\Pi$ of the flow $\phi$ on $V$. Concerning the measure $C_f$, whose support is localized at the arrows $(b,t)\in G_+^*$ with $b\in B$ a fixed point of $\phib$ and $t\in \rr_+^*$, we arrive at
$$
{C_f}\circ\Res( T^{-1} [\ln|Q|,T]) = \sum_{\substack{b\in B \\ \mathrm{fixed}}}\ \sum_{\substack{y\in M_b \\ \mathrm{fixed}}} \int_0^{\infty} \frac{\tr_s(\rho(e)(y,(b,t)))}{|(1-e^{\kappa_bt})\det(1-h_t'(y))|}\, dt
$$
where $h_t'$ is the action of the flow $\phi_t$ on the tangent space to the fiber $M_b$ at the fixed point $y$, and $\kappa_b$ is the exponent of the flow $\phib$ at the fixed point $b$. One has $(1-e^{\kappa_bt})\det(1-h_t'(y)) = \det(1-T\phi_t(y))$ with $T\phi_t(y)$ the flow on the tangent space $T_yM$. The latter is of the form $T\phi_t(y) = e^{\kappa_yt}$ with $\kappa_y$ the generator of the tangent flow at $y$. To end the computation, we remark that $\rho(e)(y,(b,t)) = c(y)^2 e(r(y),t)$ where $c\in \cinfc(M)$ is the cut-off function and $r(y)\in V$. Since $\sum_{y\in r^{-1}(v)}c(y)^2=1$ for any fixed point $v\in V$ of the flow $\phi$, one gets
$$
{C_f}\circ\Res( T^{-1} [\ln|Q|,T]) = \sum_{v\ \mathrm{fixed}} \int_0^{\infty} \frac{\tr_s(e^{j_vt})}{|\det(1-e^{\kappa_vt})|}\,e(v,t)\, dt
$$
where $j_v$ is the generator of the flow on the fiber $E_v$ and $\kappa_v$ is the generator of the flow on the tangent space $T_vV$.  \cqfd

\end{document}